\documentclass[12pt,thmsa,emstex]{article}
\usepackage{amsfonts}
\usepackage{t1enc}
\usepackage{amssymb}
\usepackage{epsfig}
\usepackage[french,english]{babel}
\usepackage{amsmath}
\usepackage{graphics}
\usepackage{layout}
\usepackage{latexsym}
\usepackage{color}
\usepackage{verbatim}
\usepackage{url}

\def\mylabel#1{\label{#1}}

\newtheorem{theorem}{Theorem}[section]
\newtheorem{lemma}[theorem]{Lemma}
\newtheorem{corollary}[theorem]{Corollary}
\newtheorem{proposition}[theorem]{Proposition}

\newtheorem{exercise}[theorem]{Exercise}

\newtheorem{remark}{Remark}
\newtheorem{example}{\bf{Example}}

\def\bit{\begin{itemize}}
\def\eit{\end{itemize}}

\reversemarginpar   

\def\bc{\begin{center}}
\def\ec{\end{center}}

\def\bthm{\begin{theorem}}
\def\ethm{\end{theorem}}

\def\bcor{\begin{corollary}}
\def\ecor{\end{corollary}}

\def\bprop{\begin{proposition}}
\def\eprop{\end{proposition}}

\def\blem{\begin{lemma}}
\def\elem{\end{lemma}}

\def\bex{\begin{example}}
\def\eex{\end{example}}

\def\bexo{\begin{exercise}}
\def\eexo{\end{exercise} }

\def\brem{\begin{remark}}
\def\erem{\end{remark}}
\def\prf{{\bf Proof: }}

\def\bdes{\begin{description}}
\def\edes{\end{description}}

\def\ita{\item[(a)]}
\def\itb{\item[(b)]}

\def\itc{\item[(c)]}
\def\itd{\item[(d)]}

\def\iti{\item[(i)]}
\def\itii{\item[(ii)]}

\def\itiii{\item[(iii)]}
\def\itiv{\item[(iv)]}

\def\itv{\item[(v)]}

\def\beq{\begin{equation}}
\def\eeq{\end{equation}}

\def\ben{\begin{enumerate}}
\def\een{\end{enumerate}}

\def\beqar{\begin{eqnarray}}
\def\eeqar{\end{eqnarray}}

\def\beqarr{\begin{eqnarray*}}
\def\eeqarr{\end{eqnarray*}}

\def\qed{\hfill $\Box$ \\[2ex]}
\def\prf{{\bf Proof: }\hspace{.1in}}

\catcode`\@=11

\newcommand{\E}{\mathcal{E}}
\newcommand{\F}{\mathcal{F}}

\newcommand{\x}{\textbf{x}}
\newcommand{\y}{\textbf{y}}

\def\Ind{{\mathbf 1}}
\def\RR{{\mathbb R}}  
\def\Rp{{\mathbb R}_+}   

\def\NN{{\mathbb N}}

\def\Pr{{\mathbb P}}

\def\rar{\rightarrow}
\def\eps{\varepsilon}

\begin{document}


\title{Lotka Volterra with randomly fluctuating environments \\or\\ "how switching between beneficial environments can make survival harder"\thanks{ This is a revised version of a paper previously entitled {\em Lotka Volterra in a fluctuating environment or "how good can be bad"} } }

\author{ Michel Bena\"{i}m\thanks{Institut de Math\'{e}matiques, Universit\'{e} de Neuch\^{a}tel, Rue Emile-Argand, Neuch\^{a}tel, Suisse-2000. (michel.benaim@unine.ch).} and Claude Lobry\thanks{EPI Modemic Inria and Universit\'{e} de Nice Sophia-Antipolis} }

\maketitle

\begin{abstract}
We consider two dimensional Lotka-Volterra systems in a fluctuating environment. Relying on recent results on stochastic persistence and piecewise deterministic Markov processes, we show that random switching between two environments that are both favorable to the same  species can lead to the extinction of this species or coexistence of the two competing species.
\end{abstract}
\paragraph{MSC:} 60J99; 34A60
\paragraph{Keywords:} Population dynamics, Persistence, Piecewise deterministic processes, Competitive Exclusion, Markov processes
{\small \tableofcontents}
\newpage
\section{Introduction}
\label{sec:intro}
In ecology, the principle of {\em competitive exclusion} formulated by Gause \cite{Gause32} in 1932  and later popularized by Hardin \cite{Hardin60}, asserts that when two species compete with each other for the same resource, the "better" competitor will eventually exclude the other. While there are numerous  evidences (based on laboratory experiences and natural observations) supporting this principle, the observed diversity of certain communities is in apparent contradiction with Gause's law. A striking example is given by the phytoplankton which demonstrate that a number of competing species  can  coexist despite very limited resources.  As a solution to this paradox,
Hutchinson \cite{Hutchinson61}
  suggested that sufficiently frequent variations of the environment can  keep  species abundances away from the equilibria predicted by competitive  exclusion.
   Since then, the idea that
   temporal fluctuations of the environment can reverse the trend of competitive  exclusion has been widely explored in the ecology literature
    (see e.g~\cite{ChessonWarner81}, 
    \cite{AbramsHoltRoth98} and \cite{Chesson2000} for an overview and much further references).

  Our goal here is to investigate rigorously this phenomenon for a two-species Lotka-Volterra model of competition under the assumption that the environment (defined by the parameters of the model)
 fluctuates randomly between two environments that are {\bf both} favorable to the same species. We will precisely
  describe -in terms of the parameters- the range of possible behaviors and explain why counterintuitive behaviors - including coexistence of the two species, or extinction of the species favored by the environments - can occur.

Throughout, we let $\RR$ (respectively $\RR_+, \RR_+^*$)  denote the set of real  (respectively non negative, positive)  numbers.


 An {\em environment} is a pair $\E = (A,B)$ defined by two matrices \beq
\label{eq:defA}
 A = \left(
           \begin{array}{cc}
             a & b\\
             c & d \\
           \end{array}
         \right),
 B = \left(
\begin{array}{c}
 \alpha \\
 \beta \\
\end{array}
\right),
 \eeq
  where $a,b,c,d,\alpha,\beta$ are positive numbers.


The {\em two-species competitive Lotka-Volterra vector field} associated to $\E$ is the map $F_{\E} : \RR^2 \mapsto \RR^2$ defined by
\beq
\label{eq:LV}
F_\E(x,y) = \left\{ \begin{array}{c}
             \alpha x(1- a x -  b y) \\
             \beta y(1- c x - d y)
           \end{array} \right. .
\eeq
Vector field $F_{\cal E}$  induces a dynamical system on $\RR_+^2$ given by the autonomous differential equation
\beq
\label{eq:LVode}
(\dot{x}, \dot{y}) = F_\E(x,y).
\eeq
Here $x$ and $y$ represent the abundances of two species (denoted the \x-species and \y-species for notational convenience)  and (\ref{eq:LVode})
describes their interaction in environment ${\cal E}.$

Environment ${\cal E}$ is said to be {\em favorable} to species $\x$ if
$$a < c \mbox{ and } b < d.$$
In other words, the {\em intraspecific  competition} within species $\x$  (measured by the parameter $a$)
is smaller than  the {\em interspecific competition} effect of species $\x$ on species $\y$ (measured by $c$) and
the  interspecific competition effect of species $\y$ on species $\x$ is smaller that the intraspecific  competition within species $\y.$

From now on, we let $\mathsf{Env}_{\x}$ denote the set of   environments favorable to species $\x.$
The  following result  easily follows from an isocline analysis (see e.g~ \cite{HS98}, Chapter 3.3). It can be viewed as a mathematical formulation of the competitive exclusion principle.
\bprop
\label{prop:basic}
 Suppose\footnote{
The case $\E \in \mathsf{Env}_\y$ is similar with $(0,\frac{1}{d})$ in place of $(\frac{1}{a},0)$. If now $c-a$ and $d-b$ have opposite signs, then  there is a unique equilibrium $S \in \Rp^* \times \Rp^*.$ If $c - a < 0, S$ is a sink whose basin of attraction is $\Rp^* \times \Rp^*.$  If $c - a > 0,  S$  is a saddle whose stable manifold $W^s(S)$ is the graph of a smooth bijective increasing function  $\Rp^* \rar \Rp^*.$  Orbits below $W^s(S)$ converge to $(\frac{1}{a},0)$ and orbit above  converge to  $(0,\frac{1}{d})$.}
 $\E = (A,B) \in \mathsf{Env}_\x.$ Then, for every $(x,y) \in \RR_+^* \times \RR_+$ the solution to (\ref{eq:LVode}) with initial condition $(x,y)$ converges to
$\displaystyle {(\frac{1}{a},0)}$ as $t \rar \infty.$
\eprop
If one now wants to take into account temporal variations of the environment, the autonomous system (\ref{eq:LVode}) should be replaced by the non-autonomous one
\beq
\label{eq:LVodenonauto}
(\dot{x},\dot{y}) = F_{{\cal E}(t)}(x,y),
\eeq
where, for each $t \geq 0,$ ${\cal E}(t)$ is the environment at time $t.$ The story began in the mid $1970's$ with the investigation of systems living in a periodic environment (typically justified  by the seasonal or daily fluctuation of certain abiotic factors such as  temperature or sunlight).
In 1974, Koch \cite{Koch74}, formalizing Hutchinson's ideas, described a plausible  mechanism - sustained by numerical simulations - explaining  how two species which could not coexist in a constant environment can coexist when subjected to an  additional periodic kill rate  (like seasonal harvesting or seasonal reduction of the population).  More precisely, this means that $F_{{\cal E}(t)}(x,y)$ writes
$$F_{{\cal E}(t)}(x,y) = F_{{\cal E}}(x,y) - (p(t)x, q(t) y)$$
 where $\E  \in \mathsf{Env}_\x$  and $p(t), q(t)$ are periodic positive rates. In 1980, Cushing \cite{Cu80} proves rigourously  that, under suitable conditions on $\E, p$ and $q,$ such a system may have a   locally attracting periodic orbit contained in the positive quadrant  $\Rp^* \times \Rp^*$.

In the same time and independently,  de Mottoni and Schiaffino \cite{MoSc81} prove the remarkable result that, when $t \rar \E(t)$ is $T$-periodic, every solution to (\ref{eq:LVodenonauto}) is asymptotic to a $T$-periodic orbit and construct an explicit example  having  a locally attracting positive periodic orbit,  while the averaged system
(the autonomous system (\ref{eq:LVode}) obtained from (\ref{eq:LVodenonauto}) by temporal averaging) is  favorable to the $\x$-species.
Papers \cite{Cu80} and \cite{MoSc81} are complementary. The first one relies on  bifurcations theory. The second makes a crucial use of the monotonicity  properties of the Poincaré map ($x,y \mapsto (x(T),y(T)$) and
 has inspired a large amount of work on competitive dynamics (see e.g~ the discussion and the references following Corollary 5.30 in \cite{HiSm05}).

Completely different is the approach proposed by Lobry, Sciandra, and Nival in \cite{Lobry94}. Based on classical ideas in system theory, this paper  considers
 the question from the point of view of what is now called a {\em switched system} and focus on the situation where $t \rar {\cal E}(t)$
 is  piecewise constant and assumes two possible values ${\cal E}_0, {\cal E}_1 \in \mathsf{Env}_\x.$
  For instance, Figure \ref{fig1} pictures two phase portraits (respectively colored in red and blue)  associated to the environments
 $${\cal E}_0 =
  \left( \left(
           \begin{array}{cc}
             1 & 1\\
             2 & 2 \\
           \end{array}
         \right),
\left( \begin{array}{c}
 10 \\
 1 \\
\end{array}
\right) \right) \mbox{ and }
{\cal E}_1 =
  \left( \left(
           \begin{array}{cc}
             0.5 & 0.5\\
             0.65 & 0.65 \\
           \end{array}
         \right),
\left( \begin{array}{c}
 1 \\
 10 \\
\end{array}
\right) \right)$$ both favorable to  species $\x.$ In accordance with Proposition \ref{prop:basic}
we see that all the red (respectively blue)  trajectories converge to the $x$-axis
 while  a switched trajectory  like the one shown on the picture moves  away from the $x$-axis toward the upper left direction.
This was exploited in \cite{Lobry94} to  shed light on some paradoxical effect  that had not been  previously  discussed in the literature: Even when  ${\cal E}(t) \in \mathsf{Env}_\x$ for all $t \geq 0$ (which is different from the assumption that the average vector field is induced by some ${\cal E} \in  \mathsf{Env}_\x$) not only coexistence of species but also extinction of
species $\x$ can occur.
\begin{figure}
  \centering
   	\includegraphics[width=10cm]{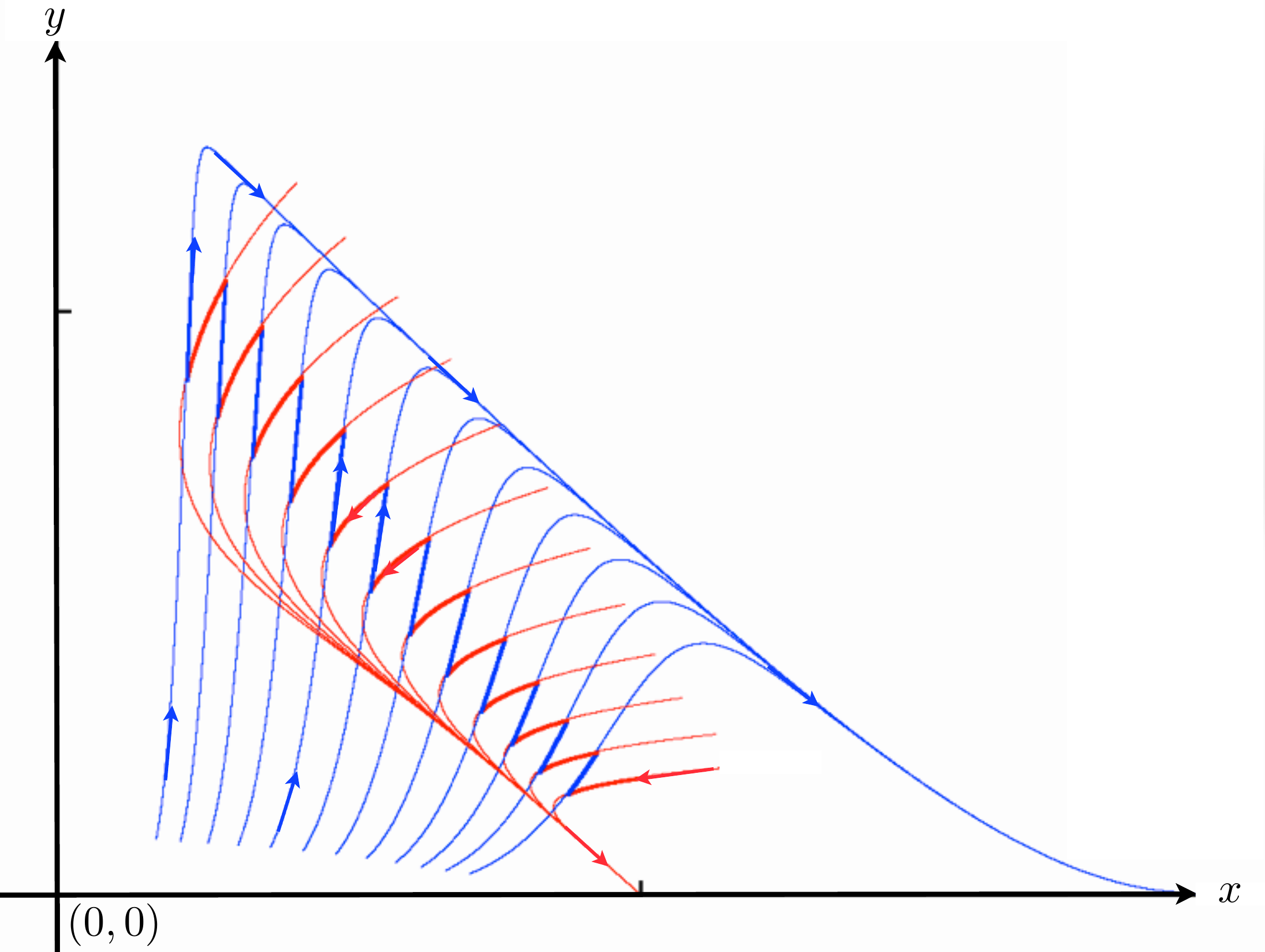}
 	\caption{An example of switched trajectory.}
 \label{fig1}
 \end{figure}

In the present paper we will pursue this line of research and investigate thoroughly the behavior of the system obtained when the environment is no longer periodic but switches randomly between ${\cal E}_0$ and ${\cal E}_1$ at jump times of a continuous time Markov chain. Our motivation is twofold: First, realistic models of environment variability should undoubtedly incorporate stochastic fluctuations. Furthermore, the  mathematical techniques involved for analyzing such a process are totally different from the deterministic ones mentioned above  and will allow to fully characterize the long term behavior of the process in terms of quantities which can be explicitly computed.

\subsection{Model, notation and presentation of main results}
\label{sec:notation}
From now on we assume given  two environments  $\E_0, \E_1 \in \mathsf{Env}_\x.$ For $i = 0, 1,$  environment $\E_i$ is defined  by (\ref{eq:defA}) with $(a_i, b_i, \ldots)$ instead of $(a,b, \ldots).$
We consider the process $\{(X_t,Y_t)\}$
defined by the differential equation
 \beq
 \label{eq:PDMP}
 (\dot{X},\dot{Y}) = F_{\E_{I_t}}(X,Y)
 \eeq
where $I_t \in \{0,1\}$ is a continuous time jump process with jump rates $\lambda_0, \lambda_1 > 0.$ That is
 $$\Pr(I_{t+s} = 1-i | I_t = i, {\cal F}_t) = \lambda_i s + o(s)$$
 where $\F_t$ is the sigma field generated by $\{I_u, u \leq t\}.$

In other words, assuming that $I_0 = i$  and $(X_0,Y_0) = (x,y),$ the process $\{(X_t,Y_t)\}$ follows the solution trajectory to $F_{\E_i}$ with initial condition $(x,y)$
for an exponentially distributed
random time, with intensity  $\lambda_i$. Then,  $\{(X_t,Y_t)\}$ follows the the solution trajectory to $F_{\E_{1-i}}$ for
another exponentially distributed random time, with intensity  $\lambda_{1-i}$ and so on.

For $\eta > 0$ small enough, the set $$K_{\eta} = \{(x,y) \in  \Rp^2 \: : \eta \leq x + y \leq 1/\eta \}$$ is positively invariant under the dynamics induced by $F_{\E_0}$ and $F_{\E_1}.$ It then attracts every solution to (\ref{eq:PDMP}) with initial condition $(x,y) \in \Rp^2 \setminus \{0,0\}.$
Fix such  $\eta >0$ and let
$$M = K_{\eta} \times \{0,1\}.$$ Set   $Z_t = (X_t,Y_t,I_t).$ Since $Z_t$ eventually lies in $M$ (whenever $(X_0,Y_0) \neq (0,0)$) we may assume without loss of generality that $Z_0 \in M$ and we see $M$ as  the
 {\em state space} of the process $\{Z_t\}_{t \geq 0}.$

The {\em extinction set of species $\y$}  is the set
$$M_0^\y = \{(x,y,i) \in M \: : y = 0\}.$$
Extinction set of species $\x,$ denoted $M_0^\x$, is defined similarly (with $x = 0$ instead of $y = 0$) and the {\em extinction set} is defined as $$M_0 = M_0^\x \cup M_0^\y.$$
The process $\{Z_t\}$ defines an homogeneous Markov process on  $M$ leaving invariant the extinction sets $M_0^\x, M_0^\y$ and the {\em interior set} $M \setminus M_0.$

It is easily seen that $\{Z_t\}$ restricted to one of the sets $M_0^\y$ or $M_0^\x$ is positively recurrent.
In order to describe its behavior on  $M \setminus M_0$ we introduce the
{\em invasion rates} of species $\y$ and $\x$
 as
\beq
\label{eq:defLambday} \Lambda_\y = \int \beta_0 (1- c_0 x)  \mu(dx,0) + \int \beta_1 (1- c_1 x)  \mu(dx,1),
\eeq
and
\beq
\label{eq:defLambday}
 \Lambda_\x  = \int \alpha_0 (1- b_0 y)  \hat{\mu}(dy,0) + \int \alpha_1 (1- b_1 y)  \hat{\mu}(dy,1).
 \eeq
 where $\mu$ (respectively $\hat{\mu}$) denotes the invariant probability measure\footnote{Here $M_0^\y$ and $M_0^\x$ are identified with  $[\eta, 1/\eta] \times \{0,1\}$ so that $\mu$ and $\hat{\mu}$  are measures on $\Rp^* \times \{0,1\}.$} of $\{Z_t\}$ on $M_0^\y$ (respectively $M_0^\x$).

Note that the quantity $\beta_i (1- c_i x)$ is the  growth rate of species $\y$ in environment  $\E_i$ when  its abundance is zero. Hence, $\Lambda_\y$ measures the long term effect of species $\x$ on the growth rate of species $\y$ when this later has low density.
 When $\Lambda_\y$ is positive (respectively negative) species $\y$ tends to increase (respectively decrease) from low density. Coexistence criteria based on the positivity of average growth rates go back to Turelli \cite{Turelli78}
   and have been used for a variety of deterministic (\cite{H81},  \cite{GH03}, \cite{S00}) and  stochastic  (\cite{ChessonEllner89}, \cite{Chesson2000}, \cite{BHW} \cite{EvansHeningSchreiber15}) models. However, these criteria are seldom expressible  in terms of the parameters of the model (average growth rates are hard to compute) and typically provide only local information on the behavior of the process near the boundary.
   Here surprisingly,
    $\Lambda_\x$ and $\Lambda_\y$ can be computed and their signs fully characterize the behavior of the process.


 Our main results can be briefly summarized as follows.
 \bdes
 \iti The invariant measures $\mu, \hat{\mu}$ and the invasion rates $\Lambda_\y$ and $\Lambda_\x$ can be explicitly computed in terms of the
 parameters ${\cal E}_i, \lambda_i, i = 0,1$ (see Section \ref{sec:prelim}).
 \itii For all $u,v \in \{+,-\}$ there are environments ${\cal E}_0, {\cal E}_1 \in \mathsf{Env}_\x$
 such that $Sign(\Lambda_\x) = u$ and $Sign(\Lambda_\y) = v.$
 Thus, in view of assertion $(iii)$ below, {\em the assumption that both environments are favorable to species $\x$
 is not sufficient do determine the outcome of the competition}.
 \itiii Let $(u,v) = (Sign(\Lambda_\x) , Sign(\Lambda_\y)).$ Assume $X_0 > 0$ and $Y_0 > 0.$ Then $(u,v)$ determines the long term behavior of $\{Z_t\}$ as follows.
 \bdes \ita  $(u,v) = (+,-) \Rightarrow  \mbox{ extinction of species \y:}$

  With probability one
  $Y_t \rar 0$ and the empirical occupation measure of $\{Z_t\}$ converges to $\mu$ (see Theorem \ref{th:good}).

 \itb  $(u,v) = (-,+) \Rightarrow  \mbox{ extinction of species \x:}$

With probability one
$X_t \rar 0$  and the empirical occupation measure of $\{Z_t\}$ converges to $\hat{\mu}$ (see Theorem \ref{th:extinc1}).
\itc $(u,v) = (-,-) \Rightarrow  \mbox{ Extinction of one species:}$

 With probability one either $X_t \rar 0$ or $Y_t \rar 0.$
  The event $\{Y_t \rar 0\}$ has positive probability. Furthermore, if the initial condition $X_0$ is sufficiently small
  or $(-,+)$ is feasible\footnote{By this, we mean that there are jump rates $\lambda_0', \lambda_1'$
such that the associated invasion rates   verify $Sign(\Lambda_\x') = -$ and $Sign(\Lambda_\y') = +$}  for ${\cal E}_0, {\cal E}_1,$ then the event $\{X_t \rar 0\}$ has positive probability (see Theorem \ref{th:extinc12}).
 \itd  $(u,v) = (+,+) \Rightarrow  \mbox{ persistence:}$

There exists a unique invariant (for $\{Z_t\}$)  probability measure $\Pi$ on $M \setminus M_0$ which is  absolutely continuous with respect to the Lebesgue measure
 $dxdy \otimes (\delta_0 + \delta_1);$ and the empirical occupation measure of $\{Z_t\}$ converges almost surely to $\Pi.$ Furthermore, for generic parameters, the law of the process converge exponentially fast to $\Pi$ in total variation. (see Theorem \ref{th:fair}).

 The density of $\Pi$ cannot be explicitly computed, still its tail behavior (Theorem \ref{th:fair}, $(ii)$) and  the topological properties of its support are well understood (see Theorem \ref{th:support}).

 \edes
\edes
The proofs rely on recent results on stochastic persistence given in \cite{Ben14}
 built upon  previous results obtained for deterministic systems in  \cite{H81, S00, GH03, HofSch04} (see also \cite{ST11}
 for a comprehensive introduction to the deterministic theory), stochastic differential equations with a small diffusion term in \cite{BHW},
 stochastic differential equations and random difference equations in \cite{SBA11, Sch12}. We also make a crucial use of some recent results on piecewise deterministic Markov
 processes obtained in \cite{bakhtin&hurt, bakhtin&hurt&matt} and \cite{BMZIHP}.

The paper is organized as follows. In Section \ref{sec:prelim} we compute $\Lambda_\x$ and $\Lambda_\y$  and derive some of their main properties.
Section \ref{sec:extinct} is devoted to the situation where one invasion rate is negative and contains the  results corresponding to the cases $(iii), (a), (b), (c)$ above. Section \ref{sec:fair} is devoted to the situation where both invasion rates are positive and contains the results corresponding to $(iii), (d)$.  Section \ref{sec:illust} presents some illustrations obtained by numerical simulation and Section \ref{sec:proof} contains the proofs of some propositions stated in section \ref{sec:prelim}.

\section{Invasion rates}
\label{sec:prelim}
As previously explained, the signs of the invasion rates will prove to be crucial for characterizing
 the long term  behavior of $\{Z_t\}.$  In this section we compute  these rates and   investigate some useful properties of the maps
  $$(\lambda_0, \lambda_1) \mapsto \Lambda_{x}(\lambda_0,\lambda_1), \Lambda_{y}(\lambda_0,\lambda_1)$$   and their zero sets.

 Set  $p_i = \frac{1}{a_i}$ and $\gamma_i = \frac{\lambda_i}{\alpha_i}.$
Here, for notational convenience,  $[p_0,p_1]$ (respectively $]p_0,p_1[$) stands for the closed (respectively open) interval with boundary points $p_0, p_1$ even when $p_1 < p_0,$ and  $M_0^\y$ is seen as a subset of   $\Rp^* \times \{0,1\}.$

The following proposition characterizes the behavior of the process on the extinction set $M_0^\y.$  The proof (given in Section \ref{sec:proof}) heavily relies on the fact that the process restricted to   $M_0^\y$, reduces to a one dimensional ODE with two possible regimes for which explicit computations are possible. It is similar to some result previously obtained  in \cite{boxma} for linear systems.
\bprop
\label{th:onM01}
The process $\{Z_t = (X_t,Y_t,I_t)\}$ restricted to $M_0^\y$  has a unique invariant probability measure $\mu$ satisfying:
\bdes
\iti
If  $p_0 = p_1 = p$ $$\mu = \delta_p \otimes \nu$$
where $\nu = \frac{\lambda_0}{\lambda_1 + \lambda_0} \delta_{1} + \frac{\lambda_1}{\lambda_1 + \lambda_0} \delta_{0}.$
\itii
If $p_0 \neq p_1$
$$\mu(dx, 1) =   h_1(x) \Ind_{[p_0, p_1]}(x) dx,$$
$$\mu(dx, 0) =   h_0(x) \Ind_{[p_0, p_1]}(x) dx $$
where
$$h_1(x) = C
\frac{ p_1 |x-p_1|^{\gamma_1-1} |p_0 - x|^{\gamma_0}}{\alpha_1 x^{1 + \gamma_0 + \gamma_1}},$$
$$h_0(x) = C
\frac{ p_0 |x-p_1|^{\gamma_1} |p_0 - x|^{\gamma_0-1}}{\alpha_0 x^{1 + \gamma_0 + \gamma_1}}$$
  and $C$ (depending on $p_1,p_0,\gamma_1,\gamma_0$) is defined by the normalization condition
  $$\int_{]p_0,p_1[} (h_1(x) + h_0(x)) dx = 1.$$
  \edes
  \eprop

For all  $x \in ]p_0,p_1[$ define
  \beq
\label{eq:deftheta}
\theta(x) =\frac {|x-p_0|^{\gamma_0 -1} |p_1 - x|^{\gamma_1 -1}}{x^{1+\gamma_0 + \gamma_1}}
\eeq and
\beq
\label{eq:defH}
P(x) = [\frac{\beta_1}{\alpha_1} (1-c_1 x) (1-a_0 x) - \frac{\beta_0}{\alpha_0} (1-c_0 x) (1-a_1x)] \frac{a_1 -a_0}{|a_1 -a_0|}.
\eeq
  Recall that the {\em invasion rate} of species $\y$ is defined (see equation (\ref{eq:defLambday})) as the growth rate of species $\y$  averaged over $\mu.$ It then follows from Proposition \ref{th:onM01} that
  \bcor
\beq
\label{eq:invas1}
\Lambda_\y = \left\{
        \begin{array}{ll}
          \frac{1}{\lambda_0 + \lambda_1} (\lambda_1 \beta_0 (1-c_0 p) + \lambda_0 \beta_1 (1-c_1 p)) & \hbox{if } p_0 = p_1 = p, \\
          \displaystyle{p_0 p_1 C \int_{]p_0,p_1[} P(x) \theta(x) dx}  & \hbox{if }  p_0 \neq p_1
        \end{array}
      \right. .
\eeq

The expression for $\Lambda_\x$ is similar. It suffices in equation (\ref{eq:invas1}) to  permute $\alpha_i$ and $\beta_i,$
and to replace  $(a_i,c_i)$ by $(d_i,b_i)$ .
\ecor
\subsection{Jointly favorable environments}
\label{sec:growthrates}
For all $0 \leq s \leq 1,$  we let  $\E_s = (A_s,B_s)$  be the environment defined by
\beq
\label{eq:averode}
s F_{\E_1} + (1-s) F_{\E_0} = F_{\E_s}
\eeq
Then, with  the notation of Section \ref{sec:notation},
$$ B_s =\left(
                                     \begin{array}{c}
                                     \alpha_s\\
                                     \beta_s\\
                                     \end{array} \right) =
                                      \left(
                                     \begin{array}{c}
                                        s \alpha_1 + (1-s) \alpha_0\\
                                       s \beta_1 + (1-s) \beta_0 \\
                                     \end{array}
                                   \right)
$$
and $$A_s =  \left(
\begin{array}{ccc}
a_s & b_s \\
\\
c_s & d_s  \\
\end{array} \right) =
\left(
\begin{array}{ccc}
 \frac{s  \alpha_1 a_1 + (1-s) \alpha_0 a_0}{\alpha_s} &
  \frac{s \alpha_1 b_1 + (1-s) \alpha_0 b_0}{\alpha_s}\\
  \\
  \frac{s  \beta_1 c_1 + (1-s) \beta_0 c_0}{\beta_s}  &
   \frac{s  \beta_1 d_1 + (1-s) \beta_0 d_0}{\beta_s}  \\
   \end{array}
  \right).$$
  Environment $\E_s$ can be understood as the environment whose dynamics (i.e the dynamics induced by  $F_{\E_s}$) is the same as the one that would result from high frequency switching giving weight $s$ to $\E_1$ and weight $(1-s)$ to $\E_0.$\footnote{More precisely,  standard averaging or mean field approximation implies that the process $\{(X_u,Y_u)\}$ with initial condition $(x,y)$ and switching rates $\lambda_0 = s t , \lambda_1 = (1-s) t$ converges in distribution, as $t \rar \infty,$ to the deterministic solution of the ODE induced by $F_{\E_s}$ and initial condition $(x,y).$}

Set
 \beq
 \label{eq:defI}
 I = \{ 0 < s < 1 : \: a_s > c_s \}
 \eeq
and
 \beq
 \label{eq:defJ}
 J = \{ 0 < s < 1 : \: b_s > d_s \}.
 \eeq

It is easily checked  that $I$ (respectively $J$) is either empty or is an open interval which closure is contained in $]0,1[.$

To get a better understanding of what $I$ and $J$ represent, observe that
\begin{itemize}
  \item If $s \in I^c \cap J^c,$  then $\E_s$ is favorable to species $x;$
   \item If  $s  \in I \cap J,$ then $\E_s$ is favorable to species $y;$
   \item If $s  \in I \cap J^c,$ then $F_{\E_s}$ has a positive sink whose basin of attraction contains the positive quadrant (stable coexistence regime);
       \item If $s   \in I^c \cap J,$ then $F_{\E_s}$  has a positive saddle whose stable manifold separates the basins of attractions of $(1/a_s,0)$ and $(0,1/d_s)$ (bi-stable regime).

\end{itemize}

 We shall say that $\E_0$ and $\E_1$ are {\em jointly favorable} to species $\x$ if for all $s \in [0,1]$ environment ${\cal E}_s$ is favorable to species $x;$ or, equivalently, $I = J = \emptyset.$   We let $\mathsf{Env}_\x^{\otimes 2} \subset \mathsf{Env}_\x \times \mathsf{Env}_\x$ denote the set of jointly favorable environments to species $x.$
\brem
\label{rem:semialg}
{\rm Set $R = \frac{\beta_0 \alpha_1}{\alpha_0 \beta_1}$
and $u = \frac{s \alpha_1}{\alpha_s}.$ Then a direct computation shows that
$ \frac{s \beta_1}{\beta_s} = \frac{u}{u (1-R) + R}.$ Thus,
$$c_s - a_s = u (c_1 \frac{1}{u(1-R) + R} - a_1) + (1-u) (c_0 \frac{R}{u(1-R) + R}) - a_0)$$
$$ = \frac{A u^2 + B u + C }{u(1-R) + R}$$
with $$A  =(a_1 - a_0) (R-1),$$
$$B  = ( 2a_0 - c_0 - a_1) R + (c_1 -a_0),$$ and $$C  = (c_0-a_0) R.$$
Then \beq
\label{eq:algI}
I \neq \emptyset \Leftrightarrow \left  \{ \begin{array}{c}
                                          A \neq 0 \\
                                          \Delta = B ^2- 4  A C > 0 \\
                                          0 < \frac{- B - \sqrt{\Delta}}{2 A} < 1
                                        \end{array} \right.
                                        \eeq
The condition for $J \neq \emptyset$ is obtained by
replacing $a_i$ by $b_i$ and $c_i$ by $d_i$ in the definitions of $A, B, C$ above, $R$ being unchanged.}
\erem
\brem
{\rm  The characterization given in  Remark \ref{rem:semialg} shows that
 $\mathsf{Env}_\x^{\otimes 2}$ is a semi algebraic subset of $\mathsf{Env}_\x \times \mathsf{Env}_\x.$
 }
\erem
The following proposition  is proved in Section \ref{sec:proof}. It provides a simple expression for  $\Lambda_\y$ in the limits of  high and low frequency switching.
\bprop
 \label{th:IJ}
The map $$\Lambda_\y : \Rp^* \times \Rp^* \mapsto \RR,$$
$$\lambda_0,\lambda_1 \mapsto \Lambda_\y(\lambda_0,\lambda_1)$$ (as defined by formulae (\ref{eq:invas1})) satisfies the following properties:
\bdes
\iti
If  $I = \emptyset,$ then for all $\lambda_0, \lambda_1$   $$\Lambda_\y(\lambda_0,\lambda_1) < 0.$$

\itii For all $s \in ]0,1[$
$$\lim_{t \rar \infty} \Lambda_\y(ts, t(1-s)) = \beta_s (1- \frac{c_s}{a_s}) \left\{
                                                                              \begin{array}{ll}
                                                                                > 0 & \hbox{ if } s \in I,\\
                                                                                 = 0 & \hbox{ if } s \in  \partial I,\\
                                                                                < 0 & \hbox{ if } s \in ]0,1[ \setminus \overline{I}
                                                                              \end{array}
                                                                            \right.
$$

 $$\lim_{t \rar 0} \Lambda_\y(ts, t(1-s)) = (1-s) \beta_0 (1- \frac{c_0}{a_0}) + s \beta_1 (1 - \frac{c_1}{a_1}) < 0;$$
 \edes
 \eprop
\brem
\label{rem:IJ}
{\rm Similarly,
\bdes
\iti
If  $J = \emptyset,$ then for all $\lambda_0, \lambda_1   \: \Lambda_\x(\lambda_0,\lambda_1) > 0$
\itii For all $s \in ]0,1[$
$$\lim_{t \rar \infty} \Lambda_\x(ts, t(1-s)) = \alpha_s (1- \frac{b_s}{d_s}) \left\{
                                                                              \begin{array}{ll}
                                                                                < 0 & \hbox{ if } s \in J,\\
                                                                                 = 0 & \hbox{ if } s \in  \partial J,\\
                                                                                > 0 & \hbox{ if } s \in ]0,1[ \setminus \overline{J}
                                                                              \end{array}
                                                                            \right.
$$

 $$\lim_{t \rar 0} \Lambda_\x(ts, t(1-s)) = (1-s) \alpha_0 (1- \frac{b_0}{d_0}) + s \alpha_1 (1 - \frac{b_1}{d_1}) > 0;$$
 \edes }
\erem
 The next result follows directly from Proposition  \ref{th:IJ} and Remark \ref{rem:IJ}.
\bcor
\label{cor:IJ}
For $u,v \in \{+,-\},$ let
 $${\cal R}_{u,v} = \{\lambda_0 > 0, \lambda_1 > 0:  \: Sign(\Lambda_\x(\lambda_0,\lambda_1)) = u, Sign(\Lambda_\y(\lambda_0,\lambda_1)) = v \}.$$ Then
 \bdes
 \iti ${\cal R}_{+-} \neq \emptyset,$
 \itii $I \cap J^c \neq \emptyset \Rightarrow {\cal R}_{+,+} \neq \emptyset,$
 \itiii $J \cap  I^c \neq \emptyset \Rightarrow {\cal R}_{-,-} \neq \emptyset,$
 \itiv $I \cap J \neq \emptyset \Rightarrow {\cal R}_{-,+} \neq \emptyset.$
 \edes
\ecor
 By using Proposition \ref{th:IJ} combined with  a beautiful argument based on second order stochastic dominance Malrieu and Zitt \cite{MalrieuZitt15} recently proved the next result. It answers a question raised in the first version of the present paper.

\begin{figure}
  \centering
 	\includegraphics[width=12cm, angle =90]{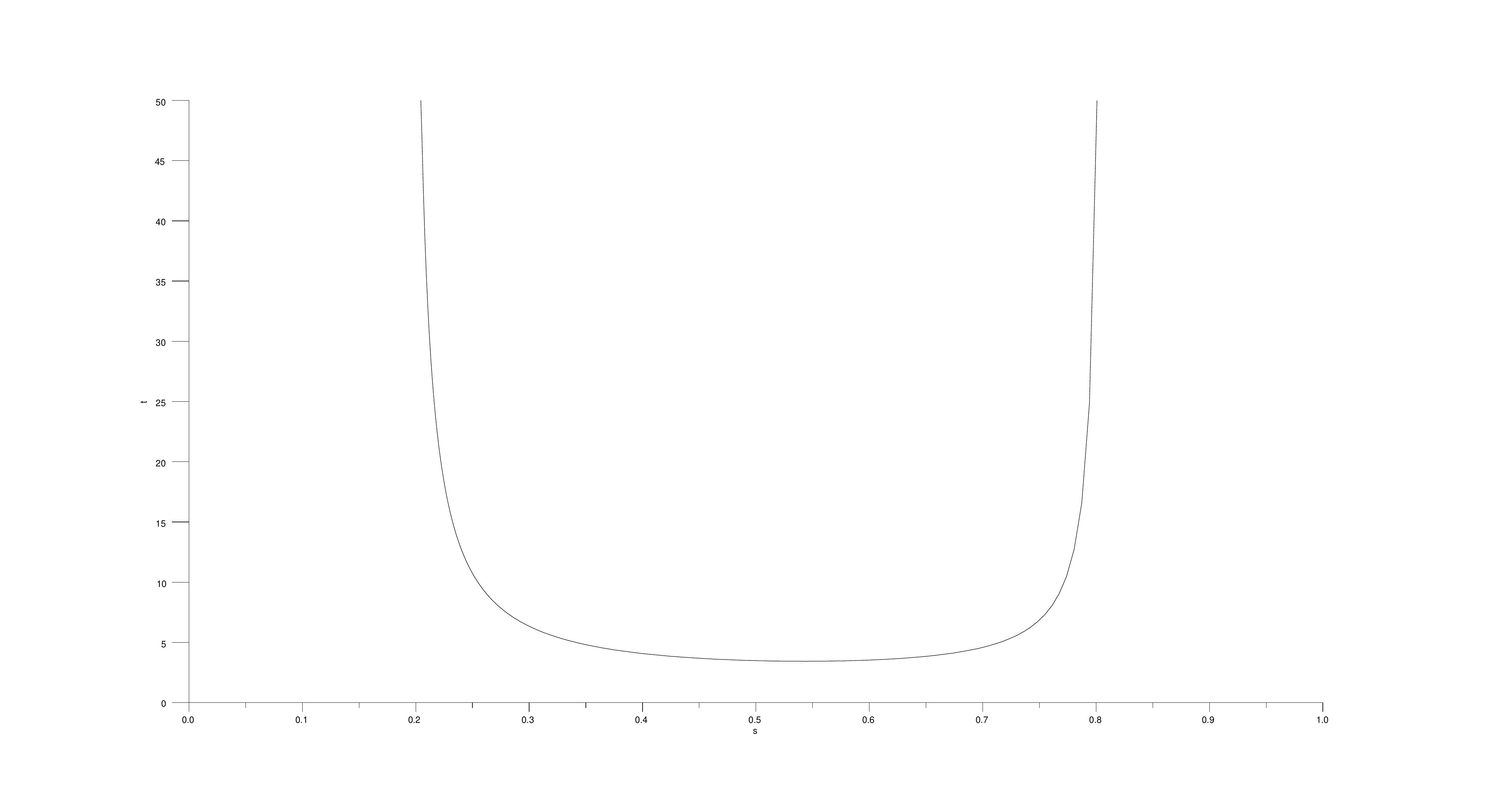}
 	\caption{Zero set of $\Lambda_{\y}(ts,(1-t)s)$  for the environments given in Section \ref{sec:illust} and $\rho = 3$ }
 \label{figlambda}
 \end{figure}

\bprop[Malrieu and Zitt, 2015]
\label{conject}
\label{th:MalZitt}
 If $I = ]s_0, s_1[ \neq \emptyset$
 the set $$\{(s,t) \in ]0,1[ \times \Rp^* \: : \Lambda_\y(ts, t(1-s)) = 0\}$$ is the graph of a continuous function $$I \mapsto \Rp^*, s \mapsto t(s)$$
 with $\lim_{s \rar s_0} t(s)= \lim_{s \rar s_1} t(s) = \infty.$
 In particular, implication $(iv)$ in Corollary \ref{cor:IJ} is an equivalence.
 \eprop
  Figure  \ref{figlambda} below represents the zero set of $s,t \mapsto \Lambda_{\y}(ts,(1-t)s)$ for the environments given in section  \ref{sec:illust} for $\rho = 3.$

 \newpage
\section{Extinction}
\label{sec:extinct}
In this section we focus on the situation where at least one invasion rate is negative and the other nonzero. If invasion rates have different signs, the species which rate is negative goes extinct and the other survives. If both are negative, one  goes extinct and the other survives.

The {\em empirical occupation} measure of the process $\{Z_t\} = \{X_t,Y_t,I_t\}$ is the (random) measure given by
$$\Pi_t = \frac{1}{t} \int_0^t \delta_{Z_s} ds.$$
Hence, for every Borel set $A \subset M, \Pi_t(A)$ is the proportion of time spent by $\{Z_s\}$ in $A$ up to time $t.$

Recall that a sequence of probability measures $\{\mu_n\}$ on a metric space $E$ (such as  $M, M_0^i$ or $\Rp^2$) is said to {\em converge weakly} to $\mu$ (another probability measure on $E$) if $\int f d\mu_n \rar \int f d\mu$ for every  bounded continuous function $f : E \mapsto \RR.$

Recall that $p_i = \frac{1}{a_i}.$
\bthm[Extinction of species $\y$]
\label{th:good}
Assume that $\Lambda_\y < 0, \Lambda_\x > 0$ and $Z_0 = z \in M \setminus M_0.$  Then, the following properties hold with probability one:
\bdes
\ita $\limsup_{t \rar \infty}  \frac{\log(Y_t)}{t} \leq \Lambda_\y,$
\itb The limit set of $\{X_t,Y_t\}$ equals $[p_0, p_1] \times \{0\},$
\itc $\{\Pi_t\}$ converges weakly to $\mu,$ where $\mu$ is the probability measure on $M_0^\y$ defined in Proposition \ref{th:onM01}
\edes
 \ethm
\brem {\rm It follows from Theorem \ref{th:good} that the marginal empirical occupation measure of $\{X_t,Y_t\}$ converges to the marginal
$$\mu(dx,0) + \mu(dx, 1) = \left\{
  \begin{array}{ll}
     \delta_{p} & \hbox{if } p_0 = p_1 = p\\
     C \theta(x) [\frac{p_1}{\alpha_1} |x-p_0| + \frac{p_0}{\alpha_0} |p_1 - x|] dx, & \hbox{if } p_0 \neq p_1
  \end{array}
\right.$$
with  $\theta$  given by (\ref{eq:deftheta}) and $C$ is a normalization constant.}
\erem
\bcor Suppose that $\E_0$ and $\E_1$ are jointly favorable to species $\x.$ Then conclusions of Theorem \ref{th:good}  hold for all positive jump rates $\lambda_0, \lambda_1.$ \ecor
\prf Follows from Theorem \ref{th:good}, Proposition \ref{th:IJ} $(i)$ and Remark \ref{rem:IJ} $(i)$. \qed
If  $\E_0$ and $\E_1$ are not jointly favorable to species $\x,$ then (by Proposition  \ref{th:IJ} and Remark \ref{rem:IJ})  there  are jump rates such that  $\Lambda_\x < 0$ or $\Lambda_{\y} > 0.$
 The following theorems tackle the situation where $\Lambda_\x < 0$.  It show that, despite the fact that environments are favorable to the same species, this species can be the one who loses the competition.
\bthm[Extinction of species $\x$]
\label{th:extinc1}
Assume that  $\Lambda_\x < 0, \Lambda_\y > 0$ and $Z_0 = z \in M \setminus M_0.$
Then, the following properties hold with probability one:
\bdes
\ita $\limsup_{t \rar \infty}  \frac{\log(X_t)}{t} \leq \Lambda_\x,$
\itb The limit set of $\{X_t,Y_t\}$ equals $\{0\} \times [\hat{p}_0, \hat{p}_1],$
\itc $\{\Pi_t\}$ converges weakly to $\hat{\mu},$ where $\hat{p}_i = \frac{1}{d_i}$ and $\hat{\mu}$ is the probability measure on $M_0^\x$ defined analogously to $\mu$ (by permuting $\alpha_i$ and $\beta_i,$
and replacing $(a_i,c_i)$ by $(d_i,b_i)$).
\edes
\ethm
\bthm[Extinction of some species]
\label{th:extinc12}
Assume that  $\Lambda_\x < 0, \Lambda_\y < 0$ and $Z_0 = z \in M \setminus M_0.$  Let  $\mathsf{Extinct}_\y$ (respectively $\mathsf{Extinct}_\x$) be the event defined by assertions $(a),(b)$ and $(c)$ in Theorem \ref{th:good} (respectively Theorem \ref{th:extinc1}). Then $$\mathbb{P}(\mathsf{Extinct}_\y) + \mathbb{P}(\mathsf{Extinct}_\x) = 1
\mbox{ and } \mathbb{P}(\mathsf{Extinct}_\y) > 0.$$
If furthermore   $z$  is sufficiently close to $M_0^\x$  or $I \cap J \neq \emptyset$ then
$$\mathbb{P}(\mathsf{Extinct}_\x) > 0.$$
\ethm

\subsection{Proofs of Theorems \ref{th:good}, \ref{th:extinc1} and \ref{th:extinc12}}
\subsubsection*{Proof of Theorem \ref{th:good}}
The strategy of the proof is the following. Assumption $\Lambda_\x > 0$ is used to show that the process eventually enter a compact set disjoint from $M_0^\x.$ Once in this compact set, it has a positive probability (independently on the starting point) to follow one of the dynamics $F_{\E_i}$ until it enters an arbitrary small neighborhood of $M_0^\y.$ Assumption $\Lambda_\y < 0$ is then used to prove that, starting from this latter neighborhood, the process  converges exponentially fast to $M_0^\y$ with positive probability. Finally, positive probability is transformed into probability one, by application of the Markov property.

Recall that $Z_t = (X_t,Y_t,I_t).$ For all $z \in M$ we let $\mathbb{P}_z$  denote the law of $\{Z_t\}_{t \geq 0}$ given that  $Z_0 = z$ and we let  $\mathbb{E}_z$ denote the corresponding expectation.

 If $E$ is one of the sets $M, M \setminus M_0, M \setminus M_0^\x$ or $M \setminus M_0^\y,$ and  $h : E \mapsto \RR$ is a measurable function which is either bounded from below or above, we let, for all  $t \geq 0$  and $z \in E,$
\beq
\label{eq:defPt}
P_t h(z) = \mathbb{E}_z(h(Z_t)).
\eeq
For $1 > \eps > 0$ sufficiently small we let
$$M_{0,\eps}^\x = \{z = (x,y,i) \in M \: : x < \eps\}$$
and
$$M_{0,\eps}^\y = \{z = (x,y,i) \in M \: : y < \eps\}$$
denote the $\eps$ neighborhoods of the extinction sets.

Let $V^\x : M \setminus M_0^\x \mapsto \RR$ and $V^\y : M \setminus M_0^\y \mapsto \RR$  be the maps defined by
$$V^\x((x,y,i)) = - \log(x) \mbox{ and } V^\y((x,y,i)) = \log(y).$$
The assumptions $\Lambda_\x > 0, \Lambda_\y < 0$ and compactness of $M_0$ imply the following Lemma:

\blem
\label{lem:lyaV12} Let $\Lambda_\x > \alpha_\x >  0$ and $ -\Lambda_\y > \alpha_\y > 0.$ Then, there exist $T > 0, \theta > 0, \eps > 0$ and
 $0 \leq \rho < 1$ such that for all $z \in M_{0,\eps}^{\textbf{h}} \setminus M_0^{\textbf{h}}, \textbf{h} \in \{\x,\y\}$
\bdes
\iti $\frac{P_T V^\textbf{h}(z) - V^\textbf{h}(z)}{T} \leq - \alpha_h,$
\itii  $P_T (e^{\theta V^\textbf{h}})(z) \leq \rho e^{\theta V^\textbf{h}}(z)$
\edes
\elem
\prf The proof can be deduced from Propositions 6.1 and 6.2  proved in a more general context in \cite{Ben14}; but for convenience and completeness we provide a simple direct proof.
We suppose $\textbf{h} = \y.$ The proof for $\textbf{h} = \x$ is identical.

$(i)$ For all $Z_0 = z \not \in M_0^{\y}$
\beq
\label{eq:basicVH}
V^\y(Z_t) - V^\y(z) = \int_0^t H(Z_s) ds
\eeq
where $$H((x,y,i) = \beta_i (1- c_i x - d_i y).$$
Thus, by taking the expectation,
$$\frac{P_T V^\y(z) - V^\y(z)}{T} = \frac{1}{T}\int_0^T P_s H(z) ds = \int H d\mu_T^z$$
where $$\mu_T^z(\cdot) = \frac{1}{T} \int_0^T P_s(z , \cdot).$$
We claim that for some $T > 0$ and $\eps > 0$  $\int H d \mu_T^z < -  \alpha_{\y} $ whenever
$z  \in  M_{0,\eps}^{\y}.$ By continuity (in $z$) it suffices to show that such a bound holds true for all $z   \in M_0^{\y}.$
By Feller continuity, compactness, and uniqueness of the invariant probability measure $\mu$  on $M_0^{\y},$ every limit point of $\{\mu_T^{z} \: :  T > 0, z  \in M_0^\y \}$ equals $\mu.$ Thus $\lim_{T \rar \infty} \int H \mu_T^z =  \int H d\mu = \Lambda_y < -\alpha_y$ uniformly in $z \in M_0^\y.$ This proves the claim and $(i).$

$(ii)$ Composing equality  (\ref{eq:basicVH})  with the map $v \mapsto e^{\theta v}$ and taking the expectation leads to
$$P_T (e^{\theta V^\y})(z) = e^{\theta V^\y(z)} e^{ l(\theta,z)}$$
where
$$l(\theta,z) = \log(\mathbb{E}_z (e^{\theta \int_0^T H(Z_s) ds})).$$
By   standard properties of the log-laplace transform, the map $\theta \mapsto l(\theta,z)$ is smooth, convex and verifies
$$l(0,z) = 0,$$
$$ \frac{\partial l}{\partial \theta}(0,z) = \mathbb{E}_z (\int_0^T H(Z_s) ds) = P_T V^\y(z) - V^\y(z)$$
and $$0 \leq \frac{\partial^2 l}{\partial \theta^2}(\theta,z) \leq \mathbb{E}_z ((\int_0^T H(Z_s) ds)^2) \leq (T \|H\|_{\infty})^2$$
where $\|H\|_{\infty} = \sup_{z \in M} |H(z)|.$
Thus, for all $z \in  M_{0,\eps}^{\y} \setminus  M_{0}^{\y}$
$$l(\theta,z) \leq T \theta (-\alpha_{\y}  + \|H\|_{\infty}^2 T \theta/2).$$
This proves $(ii),$ say for $\theta = \frac{\alpha_\y}{\|H\|_{\infty}^2 T}$ and $\rho = e^{- \frac{ \alpha_\y^2}{2 \|H\|_{\infty}^2}} .$

\qed
Define, for $\textbf{h} = \x, \y,$ the stopping times
$$\tau^{\textbf{h}, Out}_{\eps} = \min\{k \in \NN: \: Z_{kT} \in M \setminus M_{0,\eps}^\textbf{h}\}$$ and
$$\tau^{\textbf{h}, In}_{\eps} = \min\{k \in \NN: \: Z_{kT} \in  M_{0,\eps}^\textbf{h}\}.$$
{\em Step 1.} We first prove  that there exists some constant $c > 0$ such that  for all $z \in M \setminus M_0^\x$
\beq
\label{eq:claim}
\mathbb{P}_z (\tau^{\y,In}_{\eps/2} < \infty)  \geq c.
\eeq
 Set $V_k = V^\x(Z_{kT}) + k \alpha_\x T, k \in \NN.$
It follows from Lemma \ref{lem:lyaV12} $(i)$ that $\{V_{k \wedge \tau^{\x,Out}_{\eps}}\}$ is a nonnegative supermartingale. Thus, for all $z \in  M_{0,\eps}^\x \setminus M_0^\x$
$$\alpha_\x T \mathbb{E}_z (k \wedge \tau^{\x,Out}_{\eps}) \leq \mathbb{E}_z(V_{k \wedge \tau^{\x,Out}_{\eps}}) \leq V_0 = V^\x(z).$$ That is
\beq
\label{eq:tau1finite}
\mathbb{E}_z (\tau^{\x,Out}_{\eps})  \leq \frac{V^\x(z)}{\alpha_\x T} < \infty.
\eeq
Now, $(1/a_i,0)$ is  a linearly stable equilibrium for $F_{\E_i}$ whose basin of attraction contains  $\Rp^* \times \Rp$ (see Proposition \ref{prop:basic}).  Therefore, there exists $k_0 \in \NN$ such that for all $z = (x,y,i) \in M \setminus M_{0,\eps}^\x$ and $k \geq k_0$
$$\Phi^{\E_i}_{k T}(x,y) \in \{(u,v) \in \Rp  \times \Rp  \: : v < \eps/2\}.$$ Here $\Phi^{\E_i}$ stands for the flow induced by $F_{\E_i}.$ Thus, for all $z = (x,y,i) \in  M \setminus M_{0,\eps}^\x$
\beq
\label{eq:claim2}
\Pr_z( Z_{k_0 T} \in  M_{0,\eps/2}^\y) \geq \Pr(I_t = i \mbox{ for all } t \leq k_0 T | I_0 = i) = e^{-\lambda_i k_0 T} \geq c
\eeq
 where $c =  e^{-(\max{(\lambda_0, \lambda_1)} k_0 T)}.$
Combining (\ref{eq:tau1finite}) and (\ref{eq:claim2}) concludes the proof of the first step.
\\
{\em Step 2.} Let ${\cal A}$ be the event defined as
$${\cal A} = \{\limsup_{t \rar \infty} \frac{V^\y(Z_t)}{t} \leq -\alpha_\y\}.$$ We claim that there exists $c_1 > 0$ such that for all
$z \in M_{0,\eps/2}^\y$
\beq
\label{eq:claim3}
\mathbb{P}_z({\cal A}) \geq c_1.
\eeq
 Set $W_k = e^{\theta V^\y(Z_{kT})}.$ By Lemma \ref{lem:lyaV12} $(ii),$ $\{W_{k \wedge \tau^{\y,Out}_{\eps}}\}$ is a nonnegative supermartingale. Thus, for all $z \in M_{0,\eps/2}^\y$
 $$\mathbb{E}_z(W_{k \wedge \tau^{\y,Out}_{\eps}} \Ind_{\tau^{\y,Out}_{\eps} < \infty})  \leq W_0  = e^{\theta V^\y(z)} \leq (\frac{\eps}{2})^{\theta}.$$ Hence, letting $k \rar \infty$ and using  dominated convergence, leads to
 $$\eps^{\theta} \mathbb{P}_z(\tau^{\y,Out}_{\eps} < \infty) \leq \mathbb{E}_z(W_{\tau^{\y,Out}_{\eps}} \Ind_{\tau^{\y,Out}_{\eps} < \infty})  \leq W_0  \leq (\frac{\eps}{2})^{\theta}.$$
 Thus
 \beq
 \label{eq:tau2finite}
  \mathbb{P}_z(\tau^{\y,Out}_{\eps} = \infty) \geq 1- \frac{1}{2^{\theta}} = c_1 > 0.
  \eeq
 Let $M_n = \sum_{k = 1}^n (V^\y(Z_{kT}) - P_T V^\y(Z_{(k-1) T}).$ By the strong law of large numbers for  martingales applied to  $\{M_k\}$ and Lemma \ref{lem:lyaV12} (i), it follows that
 $$\limsup_{k \rar \infty} \frac{V^\y(Z_{kT})}{k T} \leq - \alpha_\y$$
 on the event $\{\tau^{\y,Out}_{\eps} = \infty \}.$ Let $C = \sup\{ \beta_i |1-c_i x - d_i y| : \: (x,y,i) \in M\}.$ It is easy to check that  $V^\y(Z_{kT + t}) - V^\y(Z_{kT})  \leq C t.$ Thus,
 $$\limsup _{t \rar \infty} \frac{V^\y(Z_{t})}{t} \leq - \alpha_\y$$ almost surely on  on the event $\{\tau^{\y,Out}_{\eps} = \infty \}.$
This later inequality, together with (\ref{eq:tau2finite}) concludes the proof of step 2.

{\em Step 3.} From (\ref{eq:claim}) and (\ref{eq:claim3}) we deduce that \beq
\label{eq:claim4}
\mathbb{P}_z({\cal A}) \geq c c_1
\eeq
for all $z \in M \setminus M_0^\x.$ Thus for all $z  \in M \setminus M_0^\x,$
 $$\Ind_{\cal A} = \lim_{t \rar \infty} \mathbb{P}_z({\cal A} | {\cal F}_t) = \lim_{t \rar \infty} \mathbb{P}_{Z_t}({\cal A}) \geq c c_1$$ where the first equality follows from  Doob's Martingale convergence theorem,  and the second from the Markov property. This concludes the proof.
\subsubsection*{Proof of Theorem \ref{th:extinc1}}
By permuting the roles of species $\x$ and $\y,$ the proof amounts to (re)proving Theorem \ref{th:good}    under the assumption that $\E_0$ and $\E_1$ are now  {\em favorable to species} $\y.$ All the arguments given in the proof of Theorem \ref{th:good} go through but for  the proof of (\ref{eq:claim2}) (where we have explicitly used the fact that  $\E_0$ and $\E_1$ are both favorable to species $\x$). In order to prove inequality (\ref{eq:claim2}) when $\E_0, \E_1 \in  \mathsf{Env}_\y$  we proceed as follows.

By Proposition \ref{th:MalZitt}  (applied after permutation of  $x$ and $y$)
the assumptions  $\E_0, \E_1 \in \mathsf{Env}_\y, \Lambda_\y < 0$ and $\Lambda_\x > 0$  imply that there exists $0 < s < 1$ such that $\E_s \in \mathsf{Env}_\x.$ Thus, there exists $k_0 \in \NN$ such that for all $z = (x,y,i) \in M \setminus M_{0,\eps}^\x$ and $k \geq k_0$
\beq
\label{eq:inside}
\Phi^{\E_s}_{k T}(x,y) \in \{(u,v) \in \Rp  \times \Rp  \: : v < \eps/2\}
\eeq
where $\Phi^{\E_s}$ stands for the flow induced by $F_{\E_s}.$
We claim that there exists $c > 0$ such that
\beq
\label{eq:claim2bis}
\Pr_z( Z_{k_0 T} \in  M_{0,\eps/2}^\y) \geq   c
\eeq
for all $z \in M \setminus M_{0,\eps}^\x.$
Suppose to the contrary that for some sequence $z_n \in M \setminus M_{0,\eps}^\x$
$$\lim_{n \rar \infty} \Pr_{z_n}( Z_{k_0 T} \in  M_{0,\eps/2}^\y) = 0.$$
By compactness of $M \setminus M_{0,\eps}^\x,$ we may assume that $z_n \rar z^* = (x^*,y^*,i^*) \in M_{0,\eps}^\x.$ Thus, by Feller continuity (Proposition 2.1 in \cite{BMZIHP}) and Portmanteau's theorem, it comes that
\beq
\label{eq:notinside}
\Pr_{z^*}( Z_{k_0 T} \in  M_{0,\eps/2}^\y) = 0.
\eeq
Now, by the support theorem (Theorem 3.4  in \cite{BMZIHP}), the deterministic orbit $\{ \Phi^{\E_s}_t(x^*,y^*): \: t \geq 0\}$ lies in the topological support of the law of $\{X_t,Y_t\}.$ This shows that (\ref{eq:notinside}) is in contradiction with (\ref{eq:inside}).
\subsubsection*{Proof of Theorem \ref{th:extinc12}}
The proof is similar to the proof of Theorem \ref{th:good}, so we only give a sketch of it.  Reasoning like in Theorem \ref{th:good}, we show that there exists $c, c_1 > 0$ such that  for all $z \in M_{0,\eps}^\textbf{h}, \mathbf{P}_z(\mathsf{Extinct}_\textbf{h}) \geq c_1$ and for all $z \in M \setminus M_{0,\eps}^\x$ $\mathbb{P}_z( \{Z_t\} \mbox{ enters } M_{0,\eps/2}^\y) \geq c.$

Thus, for all $z \in M \setminus M_0, \mathbb{P}_z(\mathsf{Extinct}_\y) + \mathbb{P}_z(\mathsf{Extinct}_\x) \geq c_1 + c c_1.$ Hence, by the Martingale argument used in the last step of the proof of Theorem  \ref{th:good}, we get that $\mathbb{P}_z(\mathsf{Extinct}_\y) + \mathbb{P}_z(\mathsf{Extinct}_\x) = 1.$ Since $(1/a_i, 0)$ is a linearly stable equilibrium for $F_{\E_i}$ whose basin contains $\Rp^* \times \Rp^*,$ $\mathbb{P}_z( \{Z_t\} \mbox{ enters } M_{0,\eps/2}^\y) > 0$ for all $z \in M \setminus M_0$ and, consequently, $\mathbb{P}_z(\mathsf{Extinct}_\y) > 0.$ If furthermore there is some $s \in I \cap J$
 $(0,1/d_s)$ is a linearly stable equilibrium for $F_{\E_s}$ whose basin contains $\Rp^* \times \Rp^*$ and, by the same argument, $\mathbb{P}_z(\mathsf{Extinct}_\y) > 0.$
\section{Persistence}
\label{sec:fair}
Here we assume that the invasion rates are positive and show that this implies a form of "stochastic coexistence".
\bthm
\label{th:fair}
 Suppose that $\Lambda_\x > 0, \Lambda_\y > 0$
 Then, there exists a unique invariant  probability measure (for the process $\{Z_t\}$) $\Pi$ on $M \setminus M_0$ i.e~ $\Pi(M \setminus M_0) = 1.$ Furthermore,
 \bdes
 \iti $\Pi$ is  absolutely continuous with respect to the Lebesgue measure $dx dy \otimes (\delta_0 + \delta_1);$
 \itii There exists $\theta > 0$ such that $$\int (\frac{1}{x^{\theta}} + \frac{1}{y^{\theta}}) d\Pi < \infty;$$
 \itiii  For every initial condition $z = (x,y,i) \in M \setminus M_0$
$$\lim_{t \rar \infty} \Pi_t = \Pi$$
weakly, with probability one.
\itiv Suppose  that $\frac{\beta_0 \alpha_1}{\alpha_0 \beta_1} \neq \frac{a_0 c_1}{a_1 c_0}$ or   $\frac{\beta_0 \alpha_1}{\alpha_0 \beta_1} \neq \frac{b_0 d_1}{b_1 d_0}.$  Then there exist constants $C, \lambda > 0$ such that  for  every Borel set $A \subset M \setminus M_0$ and every $z = (x,y,i) \in M \setminus M_0$
$$|\Pr( Z_t \in A| Z_0 = z) - \Pi(A)| \leq C (1 + \frac{1}{x^{\theta}} + \frac{1}{y^{\theta}}) e^{-\lambda t}.$$
\edes
\ethm
Theorem \ref{th:fair} has several consequences which express that, whenever the invasion rates are positive, species abundances tend to stay away from the extinction set.
 Recall that the $\eps$-boundary of the extinction set is the set
 $$M_{0,\eps} = \{z = (x,y,i) \in M :  \: \min(x,y) \leq \eps\}.$$
 Using the terminology introduced in Chesson \cite{C82}, the process is called {\em persistent in probability} if, in  the long run, densities are very likely to remain  bounded away from zero. That is
 $$\lim_{\eps \rar 0 } \limsup_{t \rar \infty} \Pr( Z_t \in M_{0,\eps} | Z_0 = z) = 0$$ for all $z \in M \setminus M_0.$
 Similarly, it is called {\em persistent almost surely} (Schreiber \cite{Sch12}) if  the fraction of time a typical population trajectory spends near the extinction set is  very small. That is $$\lim_{\eps \rar 0 } \limsup_{t \rar \infty} \Pi_t(M_{0,\eps})  = 0$$ for all $z \in M \setminus M_0.$

 By assertion $(ii)$ of Theorem \ref{th:fair} and Markov inequality
$$\Pi(M_{0,\eps})  = O(\eps^{\theta}).$$ Thus, assertion $(iii)$ implies  almost sure persistence and assertion $(iv)$ persistence in probability.

\subsection{Proof of Theorem \ref{th:fair}}
\subsubsection*{Proof of assertions $(i), (ii), (iii).$}
By Feller continuity of $\{Z_t\}$ and compactness of $M$ the sequence $\{\Pi_t\}$ is relatively compact (for the weak convergence) and every  limit point of $\{\Pi_t\}$ is an invariant probability measure (see e.g \cite{BMZIHP}, Proposition 2.4 and Lemma 2.5).

Now, the assumption that  $\Lambda_\x$ and $\Lambda_\y$ are positive, ensure that the persistence condition given in (\cite{Ben14} sections 5 and 5.2) is satisfied. Then by  the Persistence Theorem 5.1 in \cite{Ben14} (generalizing previous results in \cite{BHW} and \cite{SBA11}), every limit point of $\{\Pi_t\}$ is a probability over $M \setminus M_0$ provided $Z_0 = z \in M \setminus M_0.$
By  Lemma \ref{lem:lyaV12} $(ii)$ every such limit point satisfies the integrability condition $(ii).$

To conclude, it then suffices to show that $\{Z_t\}$ has a unique invariant probability measure   {\bf on}  $M \setminus M_0,$ $\Pi$ and that $\Pi$ is absolutely continuous with respect to $dx dy \otimes (\delta_0 + \delta_1).$

We rely on Theorem 1 in \cite{bakhtin&hurt} (see also \cite{BMZIHP}, Theorem 4.4 and the discussion following Theorem 4.5). According to this theorem, a sufficient condition ensuring both uniqueness and absolute continuity of $\Pi$ is that
\bdes
\iti There exists an {\em accessible} point $m \in \Rp^* \times \Rp^*.$
\itii The Lie algebra generated by $(F_{\E_0},F_{\E_1})$ has full rank at point $m.$
\edes
There are several equivalent formulations of accessibility (called $D$-approachability in \cite{bakhtin&hurt}). One of them, see section 3 in \cite{BMZIHP}, is  that  for every neighborhood $U$ of $m$ and every $(x,y) \in \Rp^* \times \Rp^*$ there is a solution $\eta$ to the differential inclusion $$\dot{\eta} \in conv (F_{\E_0},F_{\E_1})(\eta),$$
$$ \eta(0) = (x,y)$$ which meet $U$ (i.e $\eta(t) \in U$ for some $t > 0$). Here $conv (F_{\E_0},F_{\E_1})$ stands for the convex hull of $F_{\E_0}$ and $F_{\E_1}.$
\brem {\rm Note that here,  accessible points are  defined as points which are accessible from every point  $(x,y) \in \Rp^* \times \Rp^*$. By invariance of the boundaries, there is no point in $\Rp^* \times \Rp^*$ which is accessible from a boundary point.} \erem

 For any environment $\E,$ let  $(\Phi_t^{\E})$  denote the flow induced by $F_{\E}$ and let
 $$\gamma^+_{\E}(m) =  \{\Phi_t^{\E}(m) : \: t \geq 0\},
 \gamma^-_{\E}(m) = \{\Phi_t^{\E}(m) : \: t \leq 0\},$$
  Since $\Lambda_\y > 0, I \neq \emptyset$ by Proposition \ref{th:IJ}. Choose $s \in I.$ Then,
 point $m_s = (1/a_s,0)$ is a hyperbolic saddle equilibrium for  $F_{\E_s}$ (as defined by equation (\ref{eq:averode})) which stable manifold is the $x$-axis and which unstable manifold, denoted $W^u_{m_s}(F_{\E_s}),$ is transverse to the $x$-axis at  $m_s.$

 Now, choose an arbitrary point $m \in W^u_{m_s}(F_{\E_s}) \cap \Rp^* \times \Rp^*.$
We claim that $m$ is accessible. A standard Poincar\'e section argument shows that there exists an arc $L$ transverse to $W^u_{m_s}(F_{\E_s})$ at $m$ and a continuous maps
$P : ]p_0 -\eta_0, p_0 + \eta_0[ \times ]0,\eta_0[ \mapsto L$ such that for all
 $(x,y) \in ]p_0 -\eta_0, p_0 + \eta_0[ \times ]0,\eta_0[$ $$\gamma^+_{\E_s}(x,y) \cap L = \{P(x,y)\}$$
and $\lim_{y \rar 0} P(x,y) = m^*.$
On the other hand, for all $x > 0, y > 0,$ $$\gamma^+_{\E_0}(x,y) \cap ]p_0 -\eta_0, p_0 + \eta_0[ \times ]0,\eta_0[ \neq \emptyset$$  because  $\E_0 \in \mathsf{Env}_\x.$ This proves the claim. Now there must be some $m \in W^u_{m_s}(F_{\E_s}) \setminus \{m_s\}$ at which $F_{\E_0}(m)$ and $F_{\E_1}(m)$ span $\RR^2.$ For otherwise $W^u_{m_s}(F_{\E_s}) \setminus \{m_s\}$ would be an invariant curve for the flows $\Phi^{\E_0}$ and $\Phi^{\E_1}$ implying that $m_s = m_0 = m_1,$ hence $a_0 = a_1$ and $I = \emptyset.$
\brem
\label{rem:openGamma}
{\rm The proof above shows that the set of accessible points has nonempty interior. This will be used later in the proofs of Theorem \ref{th:fair} (iv) and \ref{th:support}.}
\erem
\subsubsection*{Proof of assertion $(iii)$.}
 The cornerstone of the proof is the following Lemma which shows that the process satisfies a certain Doeblin's condition.
 We call a point  $z_0 \in M$  a {\em Doeblin point} provided there exist a neighborhood $U_0$ of $z_0,$ positive numbers $t_0, r_0, c_0$ and a probability measure $\nu_0$ on $M$ such that for all $z \in U_0$ and $t \in [t_0, t_0 + r_0]$
\beq
\label{doeblin}
P_t(z, \cdot ) \geq c_0 \nu_0(\cdot)
\eeq
\blem
\mylabel{lem:doeblin}
\bdes
\iti There exists an accessible point  $m_0 = (x_0,y_0) \in \Rp^* \times \Rp^*,$  such that $z_0 = (m_0,0)$ (or $(m_0,1)$) is a Doeblin point.
\itii Let $\nu_0$ be the measure associated to $z_0$ given by (\ref{doeblin}). Let $K \subset M \setminus M_0$ be a compact set.  There exist positive numbers $t_K, r_K, c_K$ such that for all $z \in K$ and  $t \in [t_K, t_K + r_K]$
$$P_t(z, \cdot ) \geq c_K \nu_0(\cdot).$$
\edes
\elem
\prf  Let $\{{\cal G}_{k}, k \in \NN\}$ be the family of vector fields defined recursively by  ${\cal G}_0 = \{F_{\E_1} - F_{\E_0}\}$ and  $${\cal G}_{k +1} = {\cal G}_{k} \cup \{[G,F_{\E_0}], [G, F_{\E_1}] \: : G \in {\cal G}_{k}\}.$$
For $m \in \Rp \times \Rp,$ let ${\cal G}_{k}(m) = \{G(m) :\: G \in {\cal G}_{k}\}.$

By Theorem 4.4  in \cite{BMZIHP}, a sufficient condition ensuring that a point $z = (x,y,i) \in M$  is a Doeblin point is that ${\cal G}_{k}(m)$ spans $\RR^2$ for some $k.$ Since
${\cal G}_{1} = \{(F_{\E_1}- F_{\E_0}), [F_{\E_1}, F_{\E_0}]\}$ it then suffices to find an accessible point $m_0$ at which $(F_{\E_1}- F_{\E_0})(m_0)$ and $[F_{\E_1}, F_{\E_0}](m_0)$ are independent. Let
$$P(x,y) = \mathsf{Det} ( (F_{\E_1}- F_{\E_0})(x,y), [F_{\E_1}, F_{\E_0}](x,y)) = \sum_{\{i,j \geq 1,  3 \leq i+ j \leq 5\}} c_{ij} x^i y^j. $$ Since  the set $\Gamma$ of accessible points has non empty interior (see remark \ref{rem:openGamma}), either $P(m_0) \neq 0$ for some  $m_0 \in \Gamma$ or all the  $c_{ij}$ are identically $0.$ A direct computation (performed with the formal calculus program Macaulay2) leads to

\medskip
{\small
 \begin{tabular}{|l|l|}
  \hline
  $c_{41}$ & $- BFH + B^2 L$ \\
  $c_{32}$ & $- 2 CFH - F^2 I + BFK + 2 BCL - BEL + CFL$ \\
  $c_{23}$ & $-CEH + BEI -CFI - 2 EFI + 2CFK + C^2L$\\
  $c_{14}$ &  $-E^2 I + C E K$ \\
  $c_{31}$ & $- 2 A  F H + 2 A B L$ \\
  $c_{22}$ & $BEG - CFG - CDH - AEH + BDI - AFI - 2 DFI - BEJ +$ \\
  & $CFJ + BDK + AFK + 2 ACL + CDL - AEL$\\
  $c_{13}$ & $ - 2 DEI + 2 CDK$ \\
  $c_{21}$ & $BDG - AFG - ADH + A^2 L$\\
  $c_{12}$ & $-D^2I + CDJ - AEJ + ADK$ \\
  \hline
\end{tabular}
}
\medskip

where
$$A = \alpha_1 - \alpha_0, B = \alpha_0 a_0 - \alpha_1 a_1, C = \alpha_0 b_0 - \alpha_1 b_1,
D = \beta_1 - \beta_0, E = \beta_0 d_0 - \beta_1 d_1,$$
$$ F = \beta_0 c_0 - \beta_1 c_1, G = \alpha_0, H = -\alpha_0 a_0, I = -\alpha_0 b_0, J = \beta_0, K = -\beta_0 d_0, L = -\beta_0 c_0.$$
Under the assumption of Theorem \ref{th:fair} $a_0 \neq a_1$ so that $A$ and $B$ cannot be simultaneously null. Thus $c_{41} = c_{31} = 0$ if and only if $FH = BL$. That is
$$\frac{\beta_0 \alpha_1}{\alpha_0 \beta_1} = \frac{a_0 c_1}{a_1 c_0}.$$
Similarly $c_{14} = c_{13} = 0$ if and only if $$\frac{\beta_0 \alpha_1}{\alpha_0 \beta_1} = \frac{b_0 d_1}{b_1 d_0}.$$
This proves that the conclusion of Lemma $(i)$ holds as long as one of these two latter equalities is not satisfied.

We now prove the second assertion. Let $z_0 = (m_0,0)$ be the Doeblin point given by $(i)$, and let  $U_0, t_0, r_0, c_0, \nu_0$ be as in the definition of such a point.
Choose $p$ in the support of $\nu_0.$ Without loss of generality we can assume that $p \in K$ (for otherwise it suffices to enlarge $K$). For all $t \geq 0$ and $\delta > 0$ let
$$O(t,\delta) = \{z \in M \: : P_t(z,U_0) > \delta\}.$$ By Feller continuity and Portmanteau theorem $O(t,\delta)$ is open. Because $m_0$ is accessible, it follows from the support theorem (Theorem  3.4 in \cite{BMZIHP}) that $$M \setminus M_0 = \cup_{t \geq 0, \delta > 0} O(t,\delta).$$ Thus, by compactness, there exist $\delta > 0$ and  $0 \leq t_1 \leq \ldots \leq t_m$ such that
$$K \subset \cup_{i = 1}^m V_i$$ where $V_i = O(t_i, \delta).$
Let $l \in \{1, \ldots, m\}$ be such that $p \in V_l.$ Choose an integer $N > \frac{t_m - t_1}{r_0}$ and set $r_i = \frac{t_i - t_1}{N}.$ Then
$\tau = t_i + N (t_0 + r_i) + N t_{l}$ is independent of $i$ and for all $z \in V_i$ and $t_0 \leq t \leq t_0 + r_0$
$$P_{\tau + t}(z, \cdot) \geq \int_{U_0} P_{t_i}(z, dz_1) \int_{V_l} P_{t_0+r_i}(z_1, dz'_1) \int_{U_0} P_{t_l} (z'_1, dz_2)$$
$$\ldots \int_{V_l} P_{t_0+r_i}(z_{N}, dz'_{N}) \int_{U_0} P_{t_l} (z'_{N}, dz_{N+1}) P_t(z_{N+1}, \cdot)$$
$$\geq \delta ( c_0 \nu_0(V_l) \delta)^N c_0 \nu_0(\cdot).$$
\qed
\blem
\mylabel{lem:lyapou}
There exist positive  numbers $\theta, T, \tilde{C}$ and $0 < \rho < 1$ such that the map
$W : M \setminus M_0 \mapsto \Rp$ defined by $$W(x,y,i) = \frac{1}{x^{\theta}} + \frac{1}{y^{\theta}}$$ verifies
$$P_{nT} W \leq \rho^n W + \tilde{C}$$ for all $n \geq 1.$
\elem
\prf By Lemma \ref{lem:lyaV12} $(ii)$ there exist $0 < \rho < 1$ and $\theta, T > 0$ such that
\beq
\label{PTWleq}
P_T W \leq \rho W + \tilde{C},
\eeq
 where
 $$\tilde{C} = \sup_{z \in M \setminus M_{0, \eps} } P_T (W) - W$$ is finite by continuity of $W$ on $M \setminus M_0$ and compactness of $M \setminus M_{0, \eps}.$
So that by iterating,
$$P_{nT}  W \leq \rho^n W + \tilde{C} \sum_{k = 1}^{n-1} \rho^k  \leq  \rho^n W + \frac{\rho}{1-\rho} \tilde{C}.$$ Replacing $\tilde{C}$ by $\frac{\rho}{1-\rho} \tilde{C}$ proves the result.
\qed
To conclude the proof of  assertion $(iii)$  we then use from the classical Harris's ergodic theorem.
 Here we rely on the following version  given (an proved) in \cite{Hairer-Mattingly} :
\bthm[Harris's Theorem]
 Let ${\cal P}$ be a Markov kernel
on a measurable space $E$
assume that
\bdes
\iti
 There exists a map $W : E \mapsto [0, \infty[$ and constants $0 < \gamma < 1, \tilde{K}$ such that
${\cal P} W \leq \gamma W + \tilde{C}$
\itii
 For some $R > \frac{2 \tilde{C}}{1-\gamma}$ there exists a probability measure $\nu$ and a constant $c$ such that
${\cal P}(x,.) \geq c \nu(.)$ whenever  $W(x) \leq R.$
\edes
Then there exists a unique invariant probability $\pi$ for ${\cal P}$ and constants $C \geq 0, 0 \leq \tilde{\gamma} < 1$ such that for every bounded measurable map
$f : E \mapsto \RR$ and all $x \in E$
$$|{\cal P}^n f(x)
 - \pi f| \leq C \tilde{\gamma}^n (1 + W(x)) \|f
 \|_{\infty}.$$
 \ethm

To apply this result, set $E = M \setminus M_0, W(x,y,i) = \frac{1}{x^{\theta}} + \frac{1}{y^{\theta}}, {\cal P} = P_{nT},$ and $\gamma = \rho^n,$ where $\theta$ and $T$ are given
by Lemma \ref{lem:lyapou} and $n \in \NN^*$  remains to be chosen. Choose $R > \frac{2 \tilde{C}}{1- \rho}$ and set $K = \{z \in M \setminus M_0 : \: W(z) \leq R \}.$
By Lemma \ref{lem:doeblin}
$P_{m t}(z, \cdot) \geq c_K^m \nu_0$ for all $t \in [t_K, t_K + r_K]$ and $z \in K.$ Choose $t \in [t_K, t_K + r_K]$ such that $t/T$ is rational, and positive integers $m,n$ such that $m/n = t/T.$
Thus $P_{n T} = P_{m t} = {\cal P}$ verifies  conditions $(i),(ii)$ above of Harris's theorem.

Let $\pi$ be the invariant probability of ${\cal P}.$ For all $t \geq 0$
$\pi P_t {\cal P} = \pi {\cal P} P_t = \pi P_t$ showing that $\pi P_t$ is invariant for ${\cal P}$. Thus $\pi = \pi P_t$ so that $\pi = \Pi.$
Now for all $t > n T$ $t = k(nT) + r$ with $k \in \NN$ and $0 \leq r < nT.$ Thus
$$| P_t f (x) - \Pi f| = |{\cal P}^k P_r f - \Pi (P_r f)| \leq  C \tilde{\gamma}^k \|f - \Pi f\|_{\infty} (1 + W(x)).$$
This concludes the proof.
\subsection{The support of the invariant measure}
 We conclude this section with a theorem describing certain properties of the topological support of $\Pi.$
 Consider again the differential inclusion induced by $F_{\E_0}, F_{\E_1} :$
 \beq
 \label{eq:inclus}
 \dot{\eta}(t) \in conv(F_{\E_0}, F_{\E_1})(\eta(t))
 \eeq
 A {\em solution } to (\ref{eq:inclus}) with initial condition $(x,y)$ is an absolutely continuous function $\eta : \RR \mapsto \RR^2$ such that $\eta(0) = (x,y)$ and (\ref{eq:inclus}) holds for almost every $t \in \RR.$

 Differential inclusion (\ref{eq:inclus}) induces a {\em set valued dynamical system} $\Psi = \{\Psi_t\}$ defined by
$$\Psi_t(x,y) = \left \{\eta(t) : \eta \mbox{ is solution to } (\ref{eq:inclus}) \mbox{ with initial condition } \eta(0) = (x,y)\right \}$$
A set $A \subset \RR^2$ is called {\em strongly positively invariant} under (\ref{eq:inclus}) if $\Psi_t(A) \subset A$ for all $t \geq 0.$ It is called {\em invariant}
if for every point $(x,y) \in A$ there exists a solution $\eta$ to (\ref{eq:inclus}) with initial condition $(x,y)$ such that $\eta(\RR) \subset A.$

The omega limit set of $(x,y)$ under $\Psi$ is the set $$\omega_{\Psi}(x,y) = \bigcap_{t \geq 0} \overline{\Psi_{[t,\infty[}(x,y)}$$
As shown in (\cite{BMZIHP}, Lemma 3.9) $\omega_{\Psi}(x,y)$ is compact, connected, invariant and strongly positively invariant under $\Psi.$
 \bthm
 \label{th:support}
 Under the assumptions of Theorem \ref{th:fair}, the topological support of $\Pi$ writes  $supp(\Pi) = \Gamma \times \{0,1\}$
 where
 \bdes
 \iti $\Gamma = \omega_{\Psi}(x,y)$ for all $(x,y) \in \Rp^* \times \Rp^*.$ In particular, $\Gamma$ is compact connected strongly positively invariant and invariant under $\Psi;$
 \itii $\Gamma$ equates the closure of its interior;
   \itiii $\Gamma \cap \Rp \times \{0\} =  [p_0,p_1] \times \{0\};$
   \itiv  If $I \cap J \neq \emptyset$ then $\Gamma \cap \{0\}  \times \Rp =   \{0\} \times [\hat{p}_0,\hat{p}_1].$
   \itv $\Gamma \setminus \{0\} \times [\hat{p}_0,\hat{p}_1]$ is contractible (hence simply connected).
   \edes
 \ethm
 \prf $(i)$ Let $(m,i) \in supp(\Pi).$ By Theorem \ref{th:fair}, for every neighborhood $U$ of $m$ and every initial condition $z = (x,y,i) \in M \setminus M_0$
 $\liminf_{t \rar \infty} \Pi_t(U) > 0.$ This implies that $m \in \omega_{\Psi}(x,y)$ (compare to Proposition 3.17 (iii) in \cite{BMZIHP}).
 Conversely, let $m \in \omega_{\Psi}(x,y)$ for some $(x,y) \in \Rp^* \times \Rp^*$ and let $U$ be a neighborhood of $m.$
 Then $$\Pi(U \times \{i\}) = \int \mathbb{P}_z (Z_z \in U \times \{i\}) \Pi(dz)
 = \int \mathbb{Q}_z (U \times \{i\}) \Pi(dz)$$
 where $\mathbb{Q}_z(\cdot) = \int_0^{\infty} \mathbb{P}_z(Z_t \in \cdot) e^{-t} dt.$ Suppose $\Pi(U \times \{i\}) = 0.$ Then for some $z_0 \in supp(\Pi) \setminus M_0$ (recall that $\Pi(M_0) = 0$) $\mathbb{Q}_{z_0}(U \times \{i\}) = 0.$
Thus $\mathbb{P}_{z_0}(Z_t \in U \times \{i\}) = 0$ for almost all $t \geq 0.$ On the other hand, because $z_0 \in supp(\Pi) \subset \omega_{\Psi}(x,y)$ there exists a solution $\eta$ to (\ref{eq:inclus}) with initial condition $(x,y)$ and some  some nonempty interval $]t_1,t_2[$ such that for all $t \in ]t_1,t_2[$ $\eta(t) \in U.$ This later property combined with the support theorem (Theorem  3.4 and Lemma 3.2 in \cite{BMZIHP}) implies that $\mathbb{P}_{z_0}(Z_t \in U \times \{i\}) > 0$ for all $t \in ]t_1,t_2[.$ A contradiction.

$(ii)$ By Proposition 3.11 in \cite{BMZIHP} (or more precisely the proof of this proposition), either  $\Gamma$ has empty interior or it equates the closure of its interior.
In the proof of Theorem \ref{th:fair}, we have shown that there exists a point $m$ in the interior of $\Gamma.$

$(iii)$ Point $(p_i,0)$ lies in $\Gamma$ as a linearly stable equilibrium of $F_{\E_i}.$ By strong invariance, $[p_0,p_1] \times \{0\} \subset \Gamma.$ On the other hand, by invariance, $\Gamma \cap \Rp \times \{0\}$ is  compact and  invariant but every compact invariant set for $\Psi$ contained in $\Rp \times \{0\}$ either equals $[p_0,p_1] \times \{0\}$ or contains the origin $(0,0).$ Since the origin is an hyperbolic linearly unstable equilibrium for $F_{\E_0}$ and $F_{\E_1}$ it cannot belong to $\Gamma.$

$(iv)$ If $I \cap J \neq \emptyset$ then for any $s \in I \cap J$ $F_{\E_s}$ has a linearly stable equilibrium $m_s \in \{0\} \times [\hat{p}_0,\hat{p}_1] $ which basin of attraction contains $\Rp^* \times \Rp^*.$ Thus $m_s \in \Gamma$ proving that $\Gamma \cap \{0\} \times \Rp$ is non empty. The proof that  $\Gamma \cap \{0\} \times \Rp = \{0\} \times [\hat{p}_0,\hat{p}_1]$ is similar to the proof of assertion $(iii).$
$(v).$ Since $\Gamma$ is positively invariant under $\Phi^{\E_0}$ and $(p_0,0)$ is a linearly stable equilibrium which basin contains $\Rp^* \times \Rp,$ $\Gamma \setminus (\{0\} \times \Rp)$ is contractible to $(p_0,0).$
\qed
\section{Illustrations}
\label{sec:illust}
We present some numerical simulations illustrating the results of the preceding sections. We consider the environments
\beq
 A_0 = \left(
           \begin{array}{cc}
             1 & 1\\
             2 & 2 \\
           \end{array}
         \right),
 B_0 = \left(
\begin{array}{c}
 1 \\
 5\\
\end{array}
\right),
 \eeq
 and
\beq
 A_1 = \left(
           \begin{array}{cc}
             3 & 3\\
             4 & 4 + \rho \\
           \end{array}
         \right),
 B_1 = \left(
\begin{array}{c}
 5 \\
 1 \\
\end{array}
\right).
 \eeq
The simulations below are obtained with $$\lambda_0 = st, \lambda_1 = (1-s) t$$ for different values of
 $s \in ]0,1[, t > 0$ and $\rho \in \{0, 1, 3\}.$
Let $S(u) = \frac{u}{5 (1-u) + u}.$ Using, Remark   \ref{rem:semialg}, it is easy to check that
\bdes
\ita $I = S(]\frac{3}{4} - \frac{1}{2 \sqrt{6}}, \frac{3}{4} + \frac{1}{2 \sqrt{6}}[),$
\itb $J = I$ for $\rho = 0,$
\itc $J = S(]\frac{71}{96} - \frac{\sqrt{241}}{96}, \frac{71}{96} + \frac{\sqrt{241}}{96}[ \subset I$ for $\rho = 1,$
\itd $J = \emptyset $ for $\rho = 3.$
\edes
The phase portraits of $F_{\E_0}$ and $F_{\E_1}$ are given in Figure \ref{fig1} with $\rho = 3.$

\begin{figure}
\includegraphics[width=6cm]{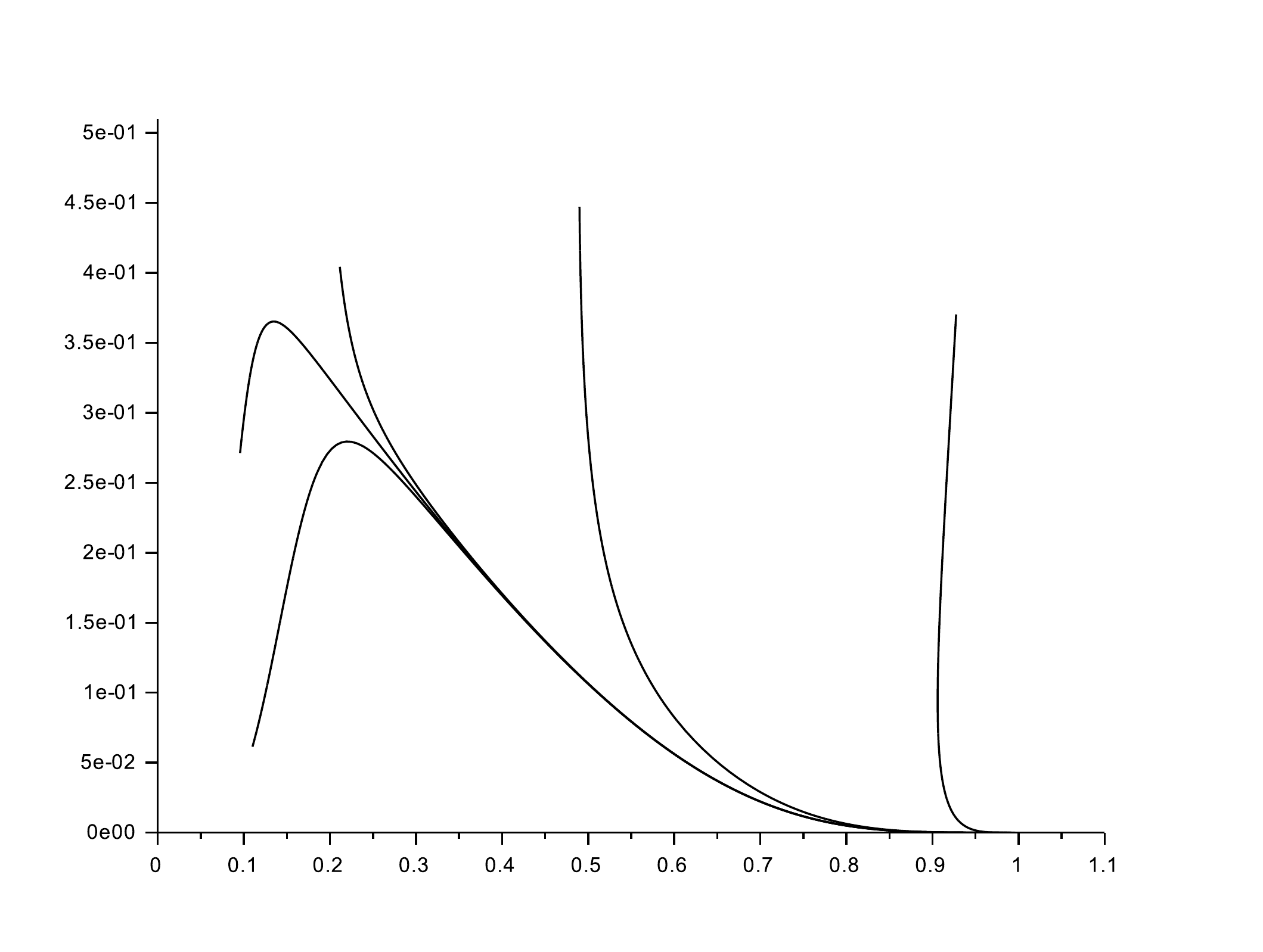}, \includegraphics[width=6cm]{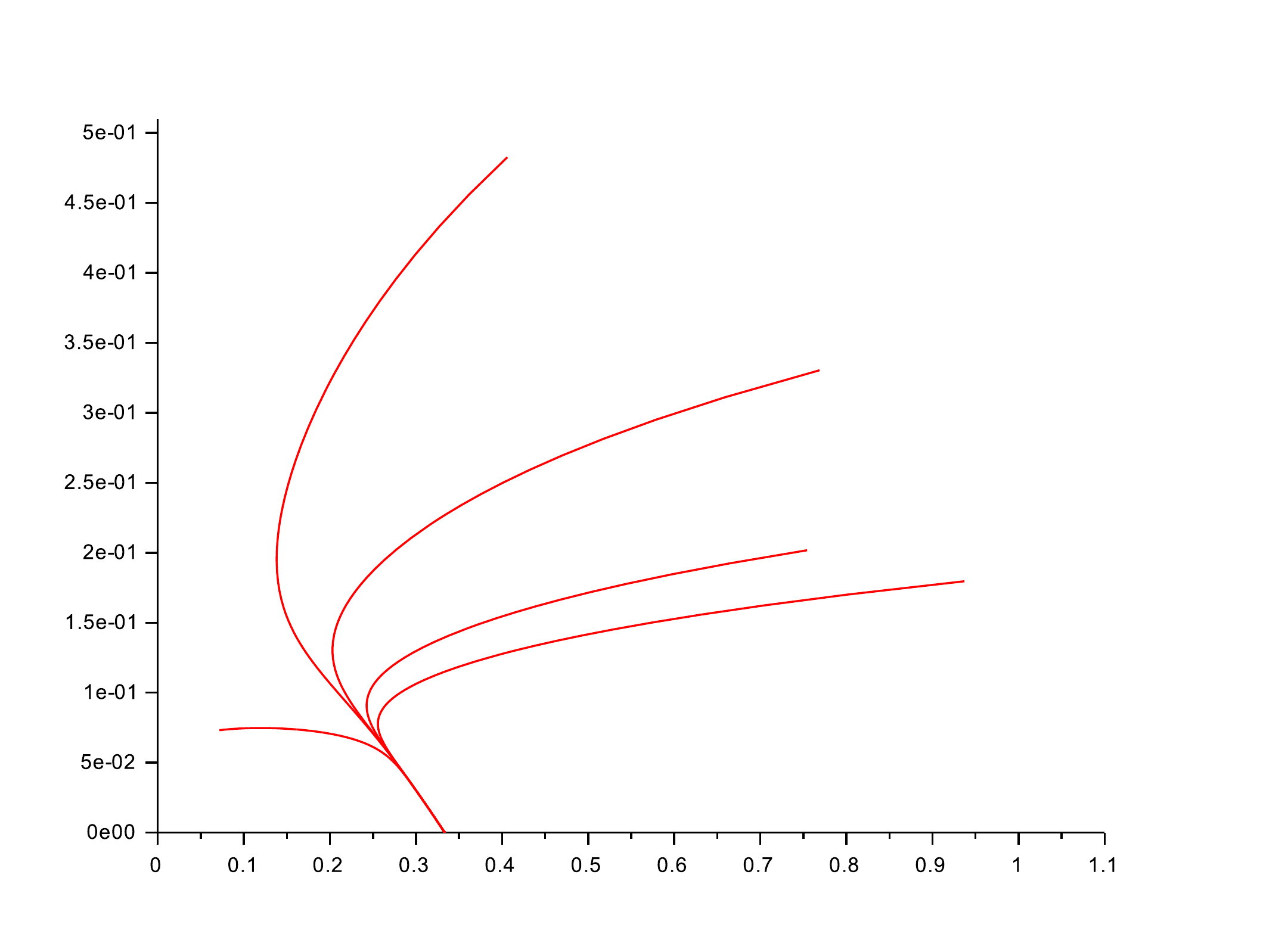}\\
  \caption{Phase portraits of $F_{\E_0}$  and $F_{\E_1}$ \label{fig1}}
\end{figure}
\begin{figure}
\centering
\includegraphics[width=10cm]{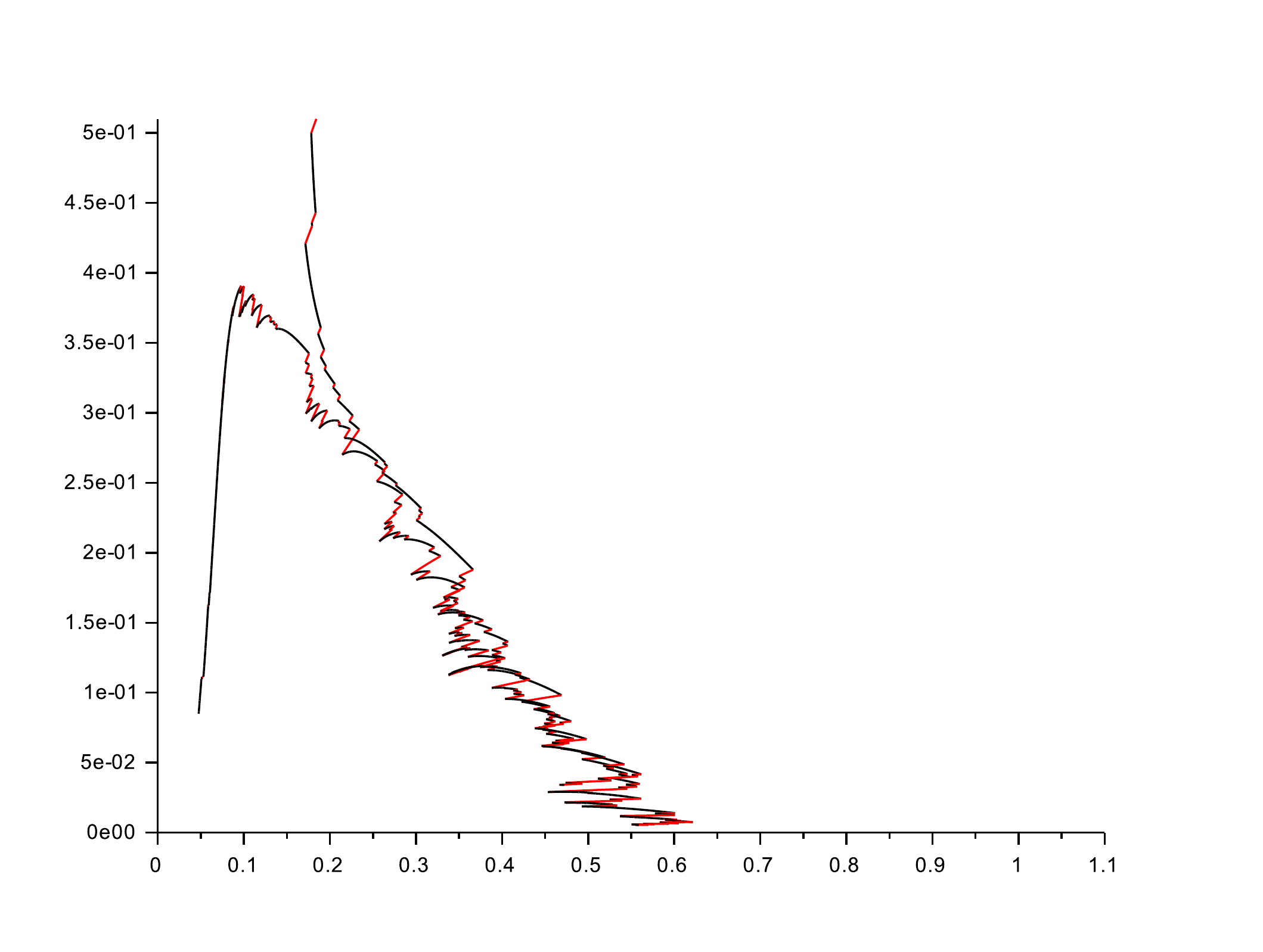}\\
  \caption{$\rho = 3, u = 0.4, t = 100$ (extinction of species $\y$) \label{fig:goodisgood}}
\end{figure}
\begin{figure}
\centering
\includegraphics[width=10cm]{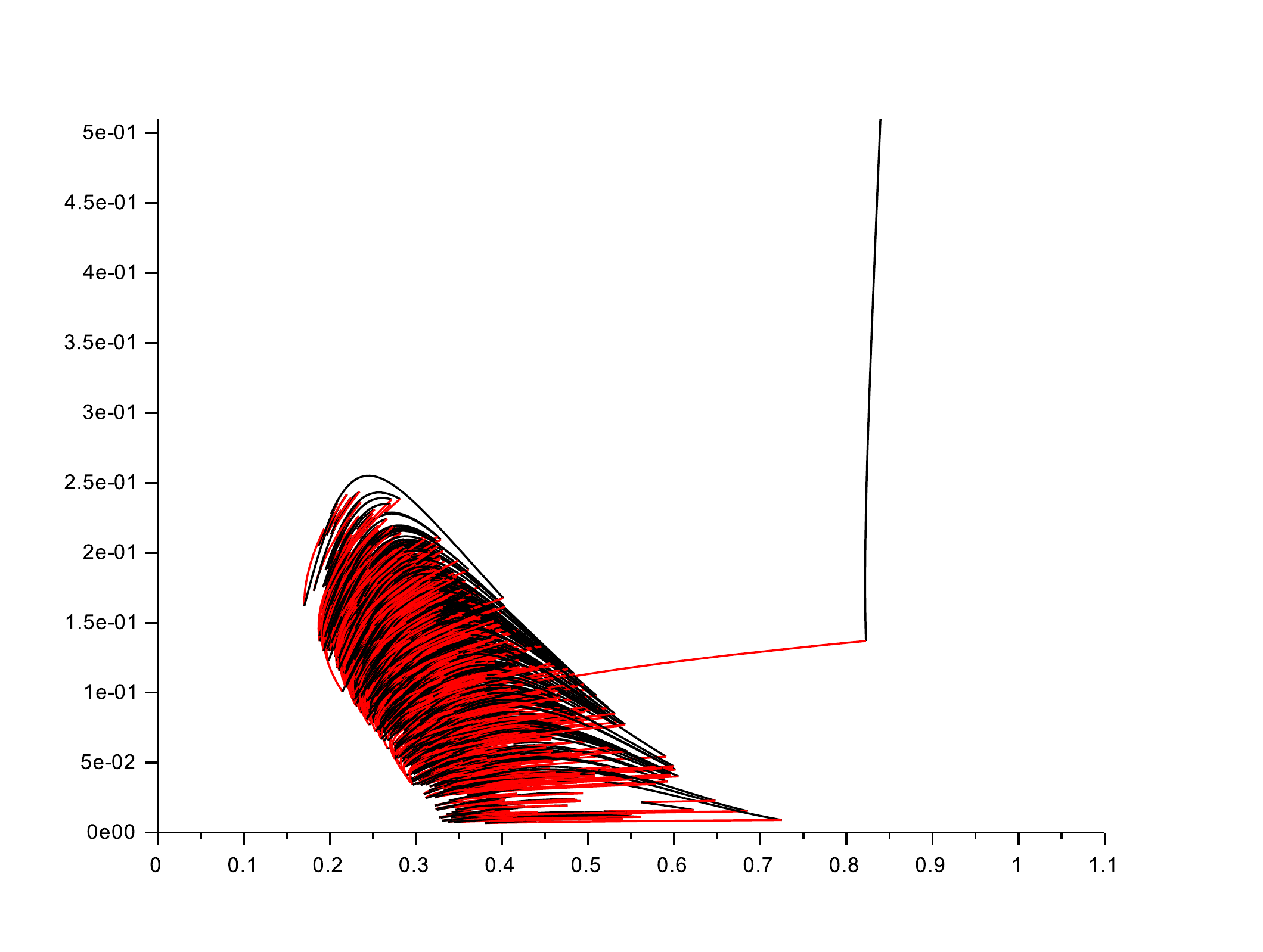}\\
  \caption{$\rho = 3, u = 0.75, t = 12$ (persistence) \label{fig:goodisfair}}
\end{figure}
\begin{figure}
\centering
  \includegraphics[width=10cm]{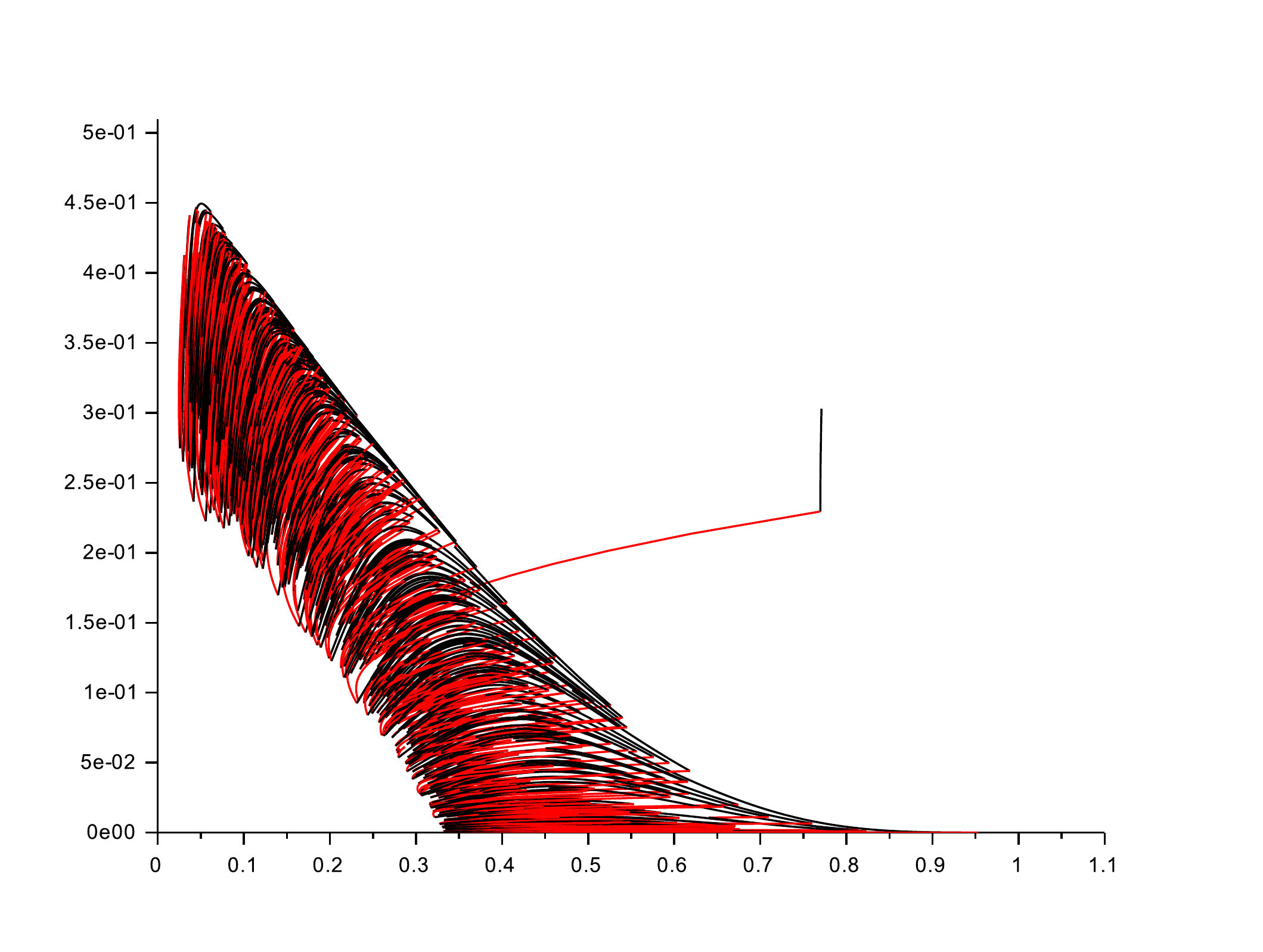}\\
  \caption{$\rho = 1, u = 0.75, t = 10$ (persistence) \label{fig:goodisfair2}}
\end{figure}
\begin{figure}
\centering
\includegraphics[width=10cm]{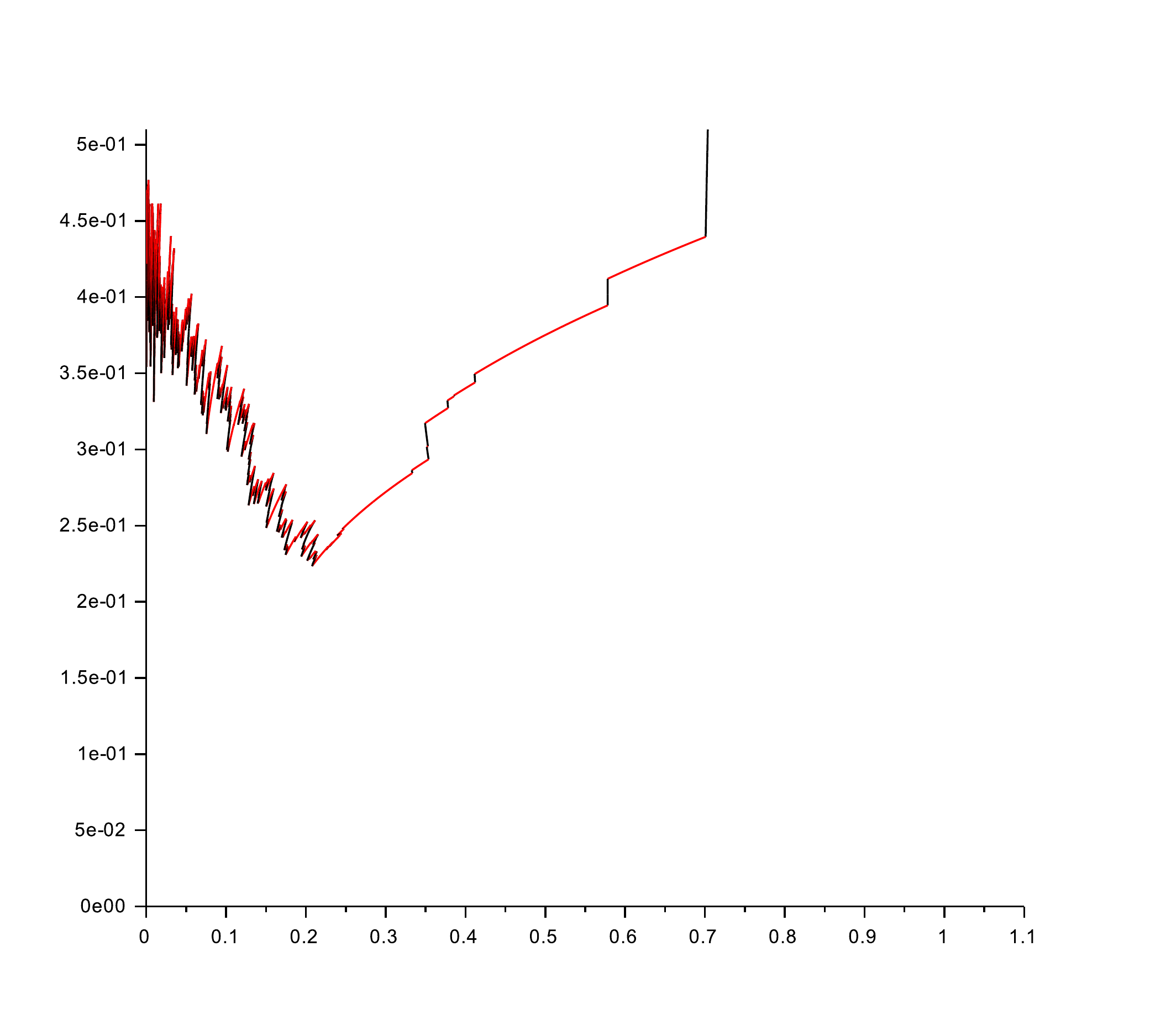}
  \caption{$\rho = 1, u = 0.75, t = 100$ (extinction of species $\x$) \label{fig:goodisbad}}
\end{figure}
\begin{figure}
\centering
\includegraphics[width=10cm]{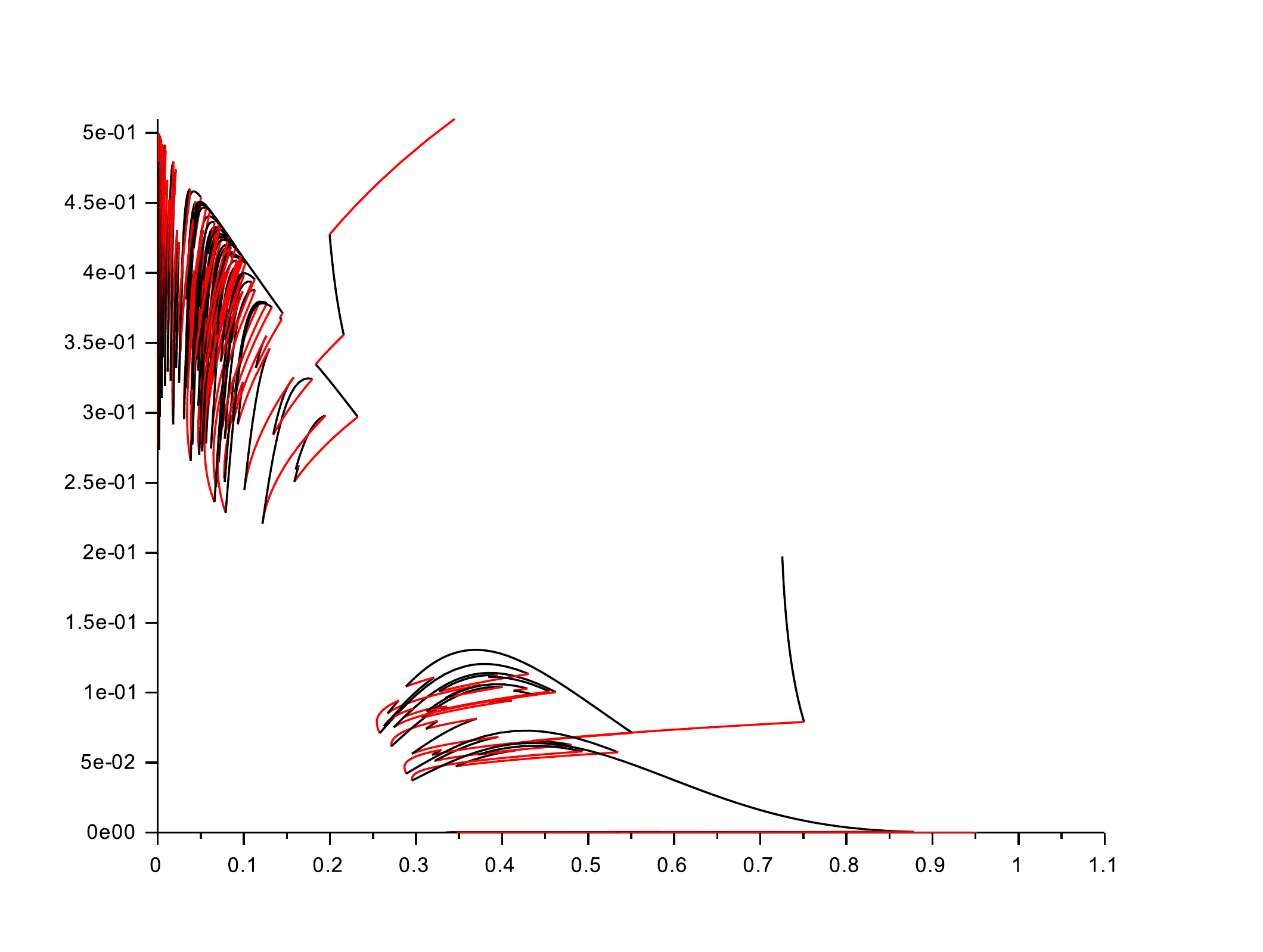}
  \caption{$\rho = 0, u = 0.75, t = 1/0.15$ (extinction of species $\x$ or $y$) \label{fig:goodisbadbad}}
\end{figure}

\paragraph{Figure  \ref{fig:goodisgood} and \ref{fig:goodisfair}} are obtained with  $\rho = 3$ (so that $J = \emptyset$).
Figure \ref{fig:goodisgood}  with $s \not \in I$ and $t$ "large" illustrates Theorems \ref{th:good} (extinction of species $\y$). Figure \ref{fig:goodisfair} with
 $s \in I$  illustrates Theorems \ref{th:fair}  and  \ref{th:support}  (persistence).
\paragraph{Figures \ref{fig:goodisfair2} and \ref{fig:goodisbad}} are obtained with $\rho = 1.$
Figure \ref{fig:goodisfair2}  with $s \in I \cap J, t = 10$ illustrates Theorems \ref{th:fair}  and   \ref{th:support} (persistence) in case $I \cap J \neq \emptyset.$  Figure \ref{fig:goodisbad} with  $s \in I \cap J$ and "large" $t$ illustrates Theorem \ref{th:extinc1}.
\paragraph{Figures \ref{fig:goodisbadbad}} is obtained with $\rho = 0$ $s \in I  \cap J$ and $t$ conveniently chosen. It illustrates Theorem \ref{th:extinc12}.
\brem
{\rm The transitions from extinction of species $\y$ to  extinction of species $\x$ when the jump rate parameter $t$ increases is reminiscent of the transition occurring with linear systems analyzed in \cite{BLMZexample} and \cite{Lawley&matt&Reed}.
}
\erem
\newpage
\section{Proofs of Propositions \ref{th:onM01} and \ref{th:IJ}}
\label{sec:proof}
\subsection{Proof of Proposition \ref{th:onM01}}
The process $\{X_t,Y_t,I_t\}$ restricted to $M_0^\y$ is defined by  $Y_t = 0$ and the one dimensional dynamics
\beq
\label{eq:logistic}
\dot{X} = \alpha_{I_t} X (1 -  a_{I_t} X)
\eeq
The invariant probability measure of the chain $(I_t)$ is given by
$$\nu = \frac{\lambda_0}{\lambda_1 + \lambda_0} \delta_{1} + \frac{\lambda_1}{\lambda_1 + \lambda_0} \delta_{0}.$$
If $a_0 = a_1 = a,  X_t \rar 1/a = p.$ Thus $(X_t,I_t)$ converges weakly to $\delta_p \otimes \nu$ and the result is proved.

Suppose now that $0 < a_0 < a_1.$

By Proposition 3.17 in \cite{BMZIHP} and Theorem 1 in \cite{bakhtin&hurt}
( or Theorem 4.4 in \cite{BMZIHP}), there exists a unique  invariant probability measure $\mu$ on $\RR_+^* \times \{0,1\}$ for $(X_t,I_t)$ which furthermore is supported by $[p_1, p_0].$
A recent result by \cite{bakhtin&hurt&matt} also proves that such a measure has a smooth density (in the $x$-variable) on $]p_1, p_0[.$

Let $\Psi : \RR \times \{0,1\} \mapsto \RR, (x,i) \mapsto \Psi(x,i)$ be smooth in the $x$ variable.
 Set $\Psi'(x,i) = \frac{\partial \Psi(x,i)}{\partial x},$ and $f_i(x) = \alpha_i x (1 - \frac{x}{p_i}).$
 The infinitesimal generator of $(x(t),I_t)$ acts on $\Psi$ as follows
 \begin{eqnarray*}
   {\cal L} \Psi(x,1) &=& \langle f_1(x), \Psi'(x , 1) \rangle + \lambda_1 (\Psi(x,0) - \Psi(x,1))  \\
   {\cal L} \Psi(x,0) &=&  \langle f_0(x), \Psi'(x , 0) \rangle + \lambda_0 (\Psi(x,1) - \Psi(x,0))
 \end{eqnarray*}
 Write $\mu(dx,1) = h_1(x) dx$ and $\mu(dx,0) = h_0(x) dx.$ Then
 $$\sum_{i = 0, 1} \int {\cal L} \Psi(x,i) h_i(x) dx = 0.$$
Choose $\Psi(x,i) = g(x) + c$ and $\Psi(x,1-i) = 0$  where $g$ is an
arbitrary compactly supported smooth function
and $c$ an arbitrary constant. Then, an easy integration by part leads to the differential equation
\beq
\label{eq:solutionh}
 \left\{
        \begin{array}{c}
               \lambda_0 h_0(x) - \lambda_1 h_1(x)  = - (f_0 h_0)'(x) \\
              \lambda_0 h_0(x)- \lambda_1 h_1(x) = (f_1 h_1)'(x)
            \end{array} \right.
 \eeq
and the condition
\beq
\label{eq:normalize}
\int_{p_1}^{p_0} \lambda_0 h_0(x) - \lambda_1 h_1(x) dx = 0.
\eeq

The maps
\beq
\label{eq:h1}
h_1(x) =  C \frac{p_1 (x-p_1)^{\gamma_1-1} (p_0 - x)^{\gamma_0}}{\alpha_1 x^{1 + \gamma_1 + \gamma_0}},\eeq
\beq
\label{eq:ho}
  h_0(x) =  C \frac{p_0 (x-p_1)^{\gamma_1} (p_0 - x)^{\gamma_0-1}}{\alpha_0 x^{1 + \gamma_1 + \gamma_0}} \eeq are solutions, where $C$ is a normalization constant  given by $$\int_{p_1}^{p_0} h_0(x) + h_1(x) dx = 1.$$ Note that $h_1$ and $h_0$ satisfy the equalities:
  $$\int_{p_1}^{p_0} h_0(x) dx = \frac{\lambda_1}{\lambda_0 + \lambda_1}$$
  $$\int_{p_1}^{p_0} h_1(x) dx = \frac{\lambda_0}{\lambda_0 + \lambda_1}.$$ This concludes the proof of Proposition \ref{th:onM01}.

\subsection{Proof of Proposition \ref{th:IJ}}
 $(i).$ We assume that  $I = \emptyset.$ If $p_0 = p_1$ then $\Lambda_\y < 0.$
Suppose  that $a_0 < a_1$ (i.e~$p_0 > p_1$) (the proof is similar for $p_0 < p_1$).
Let $p_s = \frac{1}{a_s}$ with $a_s$ being given in the definition of $A_s.$ The function $s \mapsto p_s$ maps $]0,1[$ homeomorphically onto $]p_0,p_1[$
and by definition of $\E_s$
$$s \alpha_1 (1-a_1 p_s) + (1-s) \alpha_0 (1-a_0 p_s) = 0.$$
Thus $(1-a_1 p_s) = - \frac{(1-s) \alpha_0}{s \alpha_1} (1-a_0 p_s).$ Hence
$$P(p_s)  =  \frac{(1 -a_0 p_s)}{\alpha_1 s} \beta_s (1 - c_s p_s) =  \frac{\beta_s}{\alpha_1 s} (1- a_0/a_s)(1-c_s/a_s).$$
This proves that $P(x) \leq 0$ for all  $x \in ]p_0,p_1[.$  Since $P$ is a nonzero polynomial of degree $2$, $P(x) < 0$ for all, but possibly one, points
in $]p_0,p_1[.$ Thus $\Lambda_\y < 0.$

$(ii).$ If $a_0 = a_1$ the result is obvious. Thus, we can  assume (without loss of generality) that $a_0 < a_1.$
Fix $s \in ]0,1[$ and let for all $t > 0$ $\nu_1^t$ (respectively $\nu_0^t$) be the probability measure defined as
$\nu_1^t(dx) = \frac{1}{s} h_1^t(x) \Ind_{]p_1,p_0[}(x) dx$ ($\nu_0^t(dx) = \frac{1}{1-s} h_0^t(x) \Ind_{]p_1,p_0[}(x) dx$)
where $h_1^t$ (respectively $h_0^t$) is the map defined by equation
(\ref{eq:h1}) (respectively (\ref{eq:ho})) with  $\lambda_0 = s t$ and  $\lambda_1 = (1-s) t.$
We  shall prove that
\beq
\label{eq:hiweak}
\nu_i^t \Rightarrow \delta_{p_s} \mbox{ as } t \rar \infty
\eeq
and
\beq
\label{eq:hiweak2}
\nu_i^t \Rightarrow \delta_{p_i} \mbox{ as } t \rar 0,
\eeq
where $\Rightarrow$ denotes the weak convergence.
The result to be proved follows.

Let us  prove (\ref{eq:hiweak}).
For all $x \in ]p_0,p_1[, \nu_i^t(dx) = C^t_i e^{t W(x)} [x |x - p_i|]^{-1} \Ind_{]p_1,p_0[}(x) dx$ where $C^t_i$ is a normalization constant  and
$$W(x) = \frac{s}{\alpha_0} \log(p_0 - x) + \frac{1-s}{\alpha_1} \log(x -p_1) - \frac{\alpha_s}{\alpha_0 \alpha_1} \log(x).$$
We claim  that \beq
\label{eq:laplace}
argmax_{]p_0,p_1[} W = p_s = \frac{1}{a_s}
\eeq
Indeed, set  $Q(x) = W'(x) (\alpha_0 \alpha_1 x((x-p_0)(p_1 -x)).$ It is easy to verify that
$$Q(x) = s \alpha_1(p_1-x)x - (1-s) \alpha_0 (x-p_0)x - \alpha_s (p_0 -x)(x-p_1).$$ Thus $Q(p_0) < 0, Q(p_1) > 0$ and since $Q$ is a second degree polynomial,
it suffices to show that $Q(p_s) = 0$ to conclude that $p_s$ is the global minimum of $W.$
By definition of $p_s,$ $$s \alpha_1 (1- a_1 p_s) + (1-s) \alpha_0 (1-a_0 p_s) = 0.$$ Thus
$$(1-s) \alpha_0 (p_s -p_0) = \frac{s \alpha_1 a_1 }{a_0}(p_1 - p_s).$$
Plugging this equality in the expression of $Q(p_s)$ leads to $Q(p_s) =  0.$ This proves the claim.
Now, from equation (\ref{eq:laplace}) and Laplace principle we deduce (\ref{eq:hiweak}).

We now pass to the proof of (\ref{eq:hiweak2}). It suffices to show that $\nu_i^t$ converges in probability to $p_i,$ meaning that
$\nu_i^t \{ x : \: |x-p_i| \geq \eps \} \rar 0$ as $t \rar 0.$ This easily follows from the shape of $h_i^t$ and elementary estimates.
\qed
\bibliographystyle{amsplain}
\bibliography{LotkaVolterra}

\section*{Acknowledgments}
This work was supported by the SNF grants FN 200020-149871/1 and 200021-163072/1
We thank Mireille Tissot-Daguette for her help with Scilab, Elisa Gorla for her help with Maclau2  and three anonymous referees for their useful comments and recommendations on the first version of this paper.

\end{document}